\documentclass[11pt]{article}
\usepackage{graphicx,subfigure,ifthen,color,algorithm,palatino}
\usepackage{amsmath,amsthm,amssymb,amsfonts,verbatim}
\usepackage{mathrsfs,accents,setspace}
\usepackage[autostyle]{csquotes}
\newboolean{DoubleSpaced}
\setboolean{DoubleSpaced}{false}
\def\KL{\mbox{KL}}
\def\vecof{\mbox{\rm vec}}
\def\bdeta{\boldsymbol{\eta}}
\def\bSigma{\boldsymbol{\Sigma}}
\def\mydigamma{\mbox{digamma}}
\def\dsoon{{\displaystyle\frac{1}{n}}}

\def\gothJ{\mathfrak{J}}
\def\simind{\stackrel{{\tiny \mbox{ind.}}}{\sim}}
\def\bpsi{\boldsymbol{\psi}}
\def\btheta{\boldsymbol{\theta}}
\def\by{\boldsymbol{y}}
\def\Psc{{\mathcal P}}
\def\bib{\vskip12pt\par\noindent\hangindent=1 true cm\hangafter=1}
\def\bzero{\boldsymbol{0}}
\def\smhalf{{\textstyle{\frac{1}{2}}}}
\def\real{{\mathbb R}}
\def\myand{\&\ }
\def\quarter{{\textstyle{1\over4}}}
\def\naturalNumbers{{\mathbb N}}
\def\btheta{\boldsymbol{\theta}}
\def\bI{\boldsymbol{I}}
\def\kappaREPL{\xi}
\def\ba{\boldsymbol{a}}
\def\bx{\boldsymbol{x}}
\def\bv{\boldsymbol{v}}

\def\Csc{{\mathcal C}}
\def\Dsc{{\mathcal D}}
\def\threehalves{{\textstyle{3\over2}}}
\def\oneFone{{}_1F_1}
\def\twoFone{{}_2F_1}
\def\pDens{\mathfrak{p}}
\def\pHSd{\pDens_{\mbox{\tiny{HS,$d$}}}}
\def\pHSone{\pDens_{\mbox{\tiny{HS,$1$}}}}
\def\pHStwo{\pDens_{\mbox{\tiny{HS,$2$}}}}
\def\pHSKj{\pDens_{\mbox{\tiny{HS,$K_j$}}}}
\def\sixth{\textstyle{\frac{1}{6}}}
\def\twothirds{\textstyle{\frac{2}{3}}}
\def\twofifths{\textstyle{\frac{2}{5}}}
\def\twoninths{\textstyle{\frac{2}{9}}}
\def\twotwentyfifths{\textstyle{\frac{2}{25}}}
\def\sigeps{\sigma_{\varepsilon}}

\def\bbeta{\boldsymbol{\beta}}
\def\bu{\boldsymbol{u}}
\def\gammaSUBbetaj{\gamma_{\mbox{\tiny{$\beta$}\scriptsize{$j$}}}}
\def\gammaSUBuj{\gamma_{\mbox{\scriptsize{$u$}\scriptsize{$j$}}}}
\def\dNon{d_{\bullet}}
\def\dLin{d_{\circ}}
\def\bone{\boldsymbol{1}}
\def\bX{\boldsymbol{X}}
\def\bZ{\boldsymbol{Z}}
\def\gammaGHS{\gamma_{\mbox{{\tiny GHS}}}}

\def\gammaujGHS{\gamma_{u_j,\mbox{{\tiny GHS}}}}
\def\gammabetajHS{\gamma_{\beta_j,\mbox{{\tiny HS}}}}
\def\fixOne{\tau_1^2}
\def\fixTwo{\tau_2^2}
\def\srFixOne{\tau_1}
\def\srFixTwo{\tau_2}
\definecolor{DarkOrange}{rgb}{1,0.549,0}

\newtheorem{GHSresult}{\textbf{Result}}
\begin{document}

\ifthenelse{\boolean{DoubleSpaced}}{\doublespacing}{}

\vskip5mm
\centerline{\Large\bf The Grouped Horseshoe Distribution and Its Statistical Properties}
\vskip5mm
\centerline{\normalsize\sc Virginia X. He and Matt P. Wand}
\vskip5mm
\centerline{\textit{University of Technology Sydney}}
\vskip5mm
\centerline{12th July, 2024}
\vskip5mm
\centerline{\large\bf Abstract}
\vskip2mm

The Grouped Horseshoe distribution arises from hierarchical structures 
in the recent Bayesian methodological literature aimed at selection
of groups of regression coefficients. We isolate this distribution and study 
its properties concerning Bayesian statistical inference. Most, but not all,
of the properties of the univariate Horseshoe distribution are seen
to transfer to the grouped case.


\section{Introduction}\label{sec:intro}

Since around 2010, numerous continuous distributions have been proposed 
for use as prior distributions of coefficients in Bayesian regression-type models.
Table 1 of Bai \myand Ghosh (2018) provides seven such examples, all of which
correspond to scale mixtures of Normal density functions with various polynomial-tailed
density functions. In this article we focus on one of these examples known
as the \emph{horseshoe} prior. The underlying Horseshoe distribution
(Carvalho \textit{et al.}, 2010) corresponds to the mixing distribution 
being $F_{1,1}$ for variance parameter scale mixing or Half-Cauchy 
for standard deviation scale mixing.

Most of this literature is concerned with variable selection for individual
coefficients. The \emph{grouped} extension is concerned with simultaneous 
selection of a group of variables. For example, in additive model selection
(e.g. Schiepl \textit{et al.}, 2012; He \myand Wand, 2024) 
a group of variables corresponds, typically, to a set of spline basis 
functions of a continuous predictor. Grouped variable selection is an 
attractive mechanism for deciding between the continuous predictor having 
a linear or non-linear effect. Our focus in this article is the grouped
extension of the horseshoe prior as proposed by Xu \textit{et al.} (2016).

Our first goal is determination of the underlying multivariate density function
corresponding to grouped horseshoe variable selection. This involves  
integrating out the scale mixing density function and leads to a family
of distributions, indexed by dimension, that we label the 
\emph{Grouped Horseshoe} distribution. We derive an expression
for the Grouped Horseshoe density function in terms of the
generalized exponential integral functions. As for the ordinary Horseshoe
distribution, the Grouped Horseshoe density function is shown to have
a pole at the origin.

We then investigate the grouped extensions of the various Bayesian statistical
inference properties of horseshoe priors studied by Carvalho \textit{et al.} (2010).
The score function behaviour and robustness to large signals property of
horseshoe priors, studied in Section 2 of Carvalho \textit{et al.} (2010),
is shown to extend to the grouped situation. However, the super-efficiency
property based on risk rates of convergence,
studied in Section 3.3 of Carvalho \textit{et al.} (2010), does not extend to the
grouped situation.

Our main results are laid out in Sections \ref{sec:densExplic}--\ref{sec:statProp}.
The topic of Section \ref{sec:thresh} is \emph{thresholding}, which is concerned
with practical data-based rules for deciding whether or not a coefficient 
parameter in a Bayesian model is set to zero. In this section we also
investigate use of the Grouped Horseshoe distribution for Bayesian 
generalized additive model selection as considered by the authors in He \myand Wand (2024).
Our conclusions are summarized in Section \ref{sec:conclusions}.
An online supplement provides full derivations of all results.

\subsection{Notation}

For a logical proposition $\Psc$ we let $I(\Psc)=1$ if $\Psc$ is true
and $I(\Psc)=0$ if $\Psc$ is false. The Euclidean norm of column vector $\ba$ 
is denoted by $\Vert\ba\Vert\equiv\sqrt{\ba^T\ba}$. If $\bv$ is
a random vector then $\pDens(\bv)$ denotes the density function
of $\bv$. If $f$ is a smooth function that maps $\real^d$ to $\real$ then
$\nabla_{\bx}f(\bx)$ denotes the $d\times1$ vector of partial derivatives
of $f(\bx)$ with respect to the entries of $\bx$.

\section{Density Function Explicit Form}\label{sec:densExplic}

Section 2.2 of Xu \textit{et al.} (2016) introduced the grouped horseshoe model.
The underlying distribution, which we call the Grouped Horseshoe distribution,
corresponds to setting $\sigma=\tau=G=1$ and $s_1=d$ in equation (6) of 
Xu \textit{et al.} (2016). This leads to the $d\times1$ random vector $\bx$ having 
a (standard) Grouped Horseshoe distribution if and only if
\begin{equation}
\bx|\lambda\sim N(\bzero,\lambda^2\bI_d)\quad\mbox{where}\quad
\pDens(\lambda)=\frac{2I(\lambda>0)}{\pi(1+\lambda^2)}.
\label{eq:GilbertGrape}
\end{equation}
Let $E_{\nu}$ denote the \emph{generalized exponential integral} function,
given by 
$$E_{\nu}(x)\equiv \int_1^{\infty} \exp(-xt)/t^{\nu} \,dt,\quad x,\nu\in\real$$
(e.g. 8.19.3 of Olver \textit{et al.} 2023).

\begin{GHSresult}
Let $\bx$ be a $d\times1$ random vector having a Grouped Horseshoe distribution as 
defined according to (\ref{eq:GilbertGrape}). Then the density function of $\bx$,
denoted by $\pHSd(\bx)$, is
$$\pHSd(\bx)=
\frac{\Gamma\big(\smhalf(d+1)\big)}{\sqrt{2\pi^{d+2}}}\,
\exp\big(\Vert\bx\Vert^2/2\big)
E_{(d+1)/2}\big(\Vert\bx\Vert^2/2\big)\Big/\Vert\bx\Vert^{d-1},\quad
\bx\in\real^d.$$
\label{res:explicDens}
\end{GHSresult}
\noindent
A derivation of Result \ref{res:explicDens} is given in Section \ref{sec:derivResOne}
of the supplement.

\begin{figure}[h]
\centering
\includegraphics[width=0.70\textwidth]{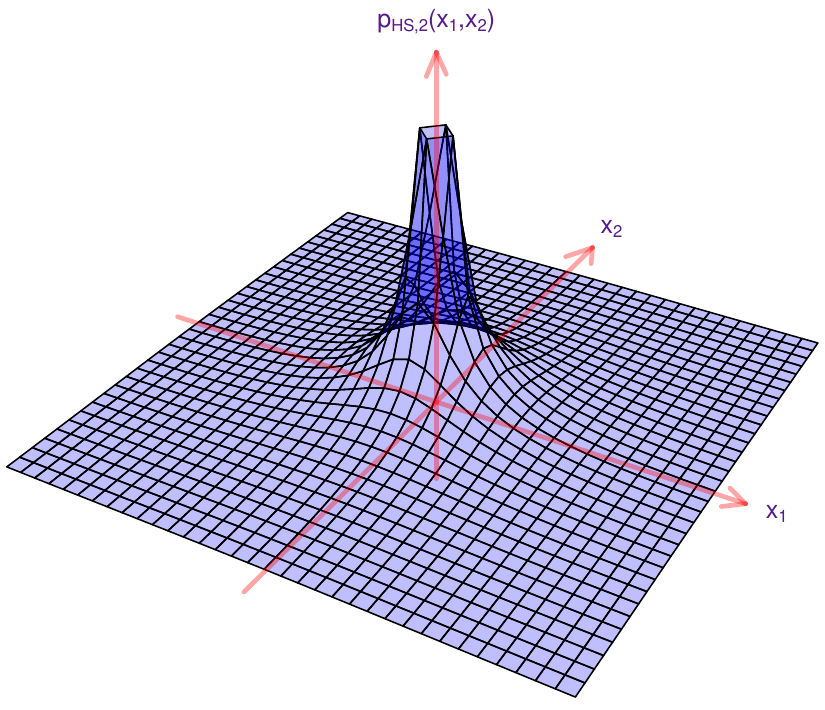}
\caption{
\textit{Perspective plot of the bivariate Grouped Horseshoe density function: $\pHStwo$.}}
\label{fig:HSbivDens}
\end{figure}

A simple consequence of Result \ref{res:explicDens} is
$$\pHSone(x)=
\frac{1}{\sqrt{2\pi^3}}\,
\exp\big(x^2/2\big)
E_1\big(x^2/2\big),\quad x\in\real.$$
which matches an expression given in the appendix of Carvalho \textit{et al.} (2010)
for the ordinary Horseshoe distribution. For $d=2$ we have 
$$\pHStwo(x_1,x_2)=
\frac{1}{2\sqrt{2\pi^3}}\,
\exp\big((x_1^2+x_2^2)/2\big)
E_{3/2}\big((x_1^2+x_2^2)/2\big)\Big/\sqrt{x_1^2+x_2^2},\quad
(x_1,x_2)\in\real^2.
$$
which is displayed in Figure \ref{fig:HSbivDens}. It is apparent 
from Figure \ref{fig:HSbivDens} that $\pHStwo$ has a pole at 
the origin. We formalise this behaviour for general $d\in\naturalNumbers$
in Section \ref{sec:pole}.

\section{Statistical Properties}\label{sec:statProp}

We now investigate various statistical properties of the Grouped Horseshoe
distribution. A particular focus is Bayesian statistical inference
where a parameter vector has a Grouped Horseshoe prior.

\subsection{Pole at the Origin Existence}\label{sec:pole}

In Carvalho \textit{et al.} (2010), $\pHSone$ is shown to have
a pole at the origin. This is shown to provide some inferential
advantages when $\pHSone$ is used as a prior density function.
Result \ref{res:poleExist}, which is derived in Section 
\ref{sec:derivResTwo} of the online supplement, shows that 
the Grouped Horseshoe density function has a pole at the origin
for any dimension.

\begin{GHSresult}
For each $d\in\naturalNumbers$,
$\displaystyle{\lim_{\bx\to\bzero}\pHSd(\bx)}=\infty$.
\label{res:poleExist}
\end{GHSresult}

\subsection{Score Function and Tail Robustness}

Consider the model
\begin{equation}
\by|\btheta\sim N(\btheta,\bI_d)\quad 
\mbox{with prior}\quad
\pDens(\btheta)=\pHSd\big(\btheta/\tau\big)\big/\tau^d
\label{eq:robMod}
\end{equation}
for some $\tau>0$ that is fixed and known. In their Section 2, 
Carvalho \textit{et al.} (2010) consider the $d=1$ version of 
(\ref{eq:robMod}) and prove that the score function 
\begin{equation}
\frac{d\log\{\pDens(y)\}}{dy}\quad\mbox{converges to $0$
as $|y|\to\infty$}.
\label{eq:ChocolateMonster}
\end{equation}
As argued there, (\ref{eq:ChocolateMonster}) implies a type of robustness
to large signals which Carvalho \textit{et al.} (2010) refer to as 
\emph{tail robustness}. 
Result \ref{res:robustRes}, which is proven in Section 
\ref{sec:derivResThree}, shows that grouped horseshoe
priors also possesses this property.

\begin{GHSresult}
For model (\ref{eq:robMod}), the tail behaviour of the score function is given by
$$\nabla_{\by}\log\{\pDens(\by)\}\sim -\frac{(d+1)\by}{\Vert\by\Vert^2}\quad\mbox{for}\quad\Vert\by\Vert\gg1.$$
Consequently
$$\lim_{\Vert\by\Vert\to\infty}\nabla_{\by}\log\{\pDens(\by)\}=\bzero
\quad\mbox{and}\quad
E\big\Vert\by-E(\btheta|\by)\big\Vert\le b_{\tau}$$
for some $b_{\tau}<\infty$ that depends on $\tau$.
\label{res:robustRes}
\end{GHSresult}

\subsection{Risk Convergence Rates}

Consider the model
\begin{equation}
\by_1,\ldots,\by_n|\btheta\quad\mbox{independently distributed as}\  
N(\btheta,\sigma^2\bI_d)
\ \mbox{with prior}\ 
\pDens(\btheta)=\pHSd\big(\btheta\big).
\label{eq:CatsMeow}
\end{equation}
Suppose that the true sampling distribution of the $\by_i$ is $N(\btheta^0,\sigma^2\bI_d)$.
In the $d=1$ case, Theorem 4 of of Carvalho \textit{et al.} (2010) states rates of convergence
results for the so-called Ces\`aro-average risk of the Bayes estimator of $\theta$, which 
they denote by $R_n$. The rates differ depending on whether $\theta^0=0$ or $\theta^0\ne0$
where $\theta^0$ is the value of $\theta$ according to the sampling distribution of the 
$y_i$. The horseshoe prior is shown to lead to a super-efficient risk rate when $\theta^0=0$.

We now provide a Grouped Horseshoe distribution extension of Theorem 4 of 
Carvalho \textit{et al.} (2010). The risk quantity $R_n$ has a definition
analogous to that given in Section 3.3 of Carvalho \textit{et al.} (2010)
for the $d$-variate extension of the set-up treated there.

\begin{GHSresult}
Consider model (\ref{eq:CatsMeow}) and suppose that the $\by_i$ have sampling
distribution $N(\btheta^0,\sigma^2\bI_d)$. Let $R_n$ be 
Ces\`aro-average risk of the Bayes estimator of $\btheta$.
When $\btheta^0=\bzero$ we have
$$
R_n\le\left\{
\begin{array}{lll}
\displaystyle{\frac{\log(n)}{2n}-\frac{\log\{\log(n)\}}{n}+O\left(\frac{1}{n}\right)}&\quad\mbox{if $d=1$},\\[2ex]
\displaystyle{\frac{\log(n)}{2n}+O\left(\frac{1}{n}\right)}&\quad\mbox{if $d\ge2$.}
\end{array}
\right.
$$
When $\btheta^0\ne\bzero$ we have $R_n\le d\log(n)/(2n)+O(1/n)$
for all $d\in\naturalNumbers$.
\label{res:riskRes}
\end{GHSresult}

Section \ref{sec:derivResFour} provides the full derivational details of Result \ref{res:riskRes}.
For $d=1$ and $\btheta^0=\bzero$, super-efficiency corresponds to presence of the $-\log\{\log(n)\}/n$
term in the upper bound on $R_n$. Result 4 shows that this term only arises in the $d=1$ case. 
The Bayes estimator is not super-efficient for $d\ge2$.

\section{Thresholding}\label{sec:thresh}

Consider a Bayesian model that contains specifications of the form
\begin{equation}
\btheta|\sigma_{\theta}
\quad\mbox{has density function}\quad
\pHSd(\btheta/\sigma_{\theta})\big/\sigma_{\theta}^d
\quad\mbox{where $\btheta$ is $d\times1$}.
\label{eq:bloodOrange}
\end{equation}
From results in Section \ref{sec:densExplic}, specification
(\ref{eq:bloodOrange}) is equivalent to 
\begin{equation}
\btheta|\sigma_{\theta},\lambda \sim N(\bzero,\sigma_{\theta}^2\lambda^2\bI_d),
\quad \pDens(\lambda)=\frac{2I(\lambda>0)}{\pi(1+\lambda^2)}
\label{eq:wetPubHol}
\end{equation}
The introduction of the auxiliary variable $\lambda$ is important
for the upcoming approach to thresholding.
In the scalar case, Carvalho \textit{et al.} (2010) develop 
a thresholding approach for deciding between 
$$\theta=0\quad\mbox{and}\quad\theta\ne 0,\quad \theta\in\real.$$
In this section we describe and evaluate the extension of their
approach to deciding between
$$\btheta=\bzero\quad\mbox{and}\quad\btheta\ne\bzero,\quad\btheta\in\real^d.$$

The Carvalho \textit{et al.} (2010) approach involves
the following result concerning a simple \enquote{side} model:

\begin{GHSresult}
For the Bayesian model
\begin{equation}
\by|\bpsi\sim N(\bpsi,\fixOne\bI_d),
\quad
\bpsi|\lambda\sim N\big(\bzero,\lambda^2\fixTwo\bI_d\big),
\quad
\pDens(\lambda)=\frac{2I(\lambda>0)}{\pi(1+\lambda^2)},
\quad\srFixOne,\srFixTwo>0\ \mbox{fixed}
\label{eq:MaineIsland}
\end{equation}
the posterior mean of $\bpsi$ is
$$E(\bpsi|\by)=E\left(\frac{\lambda^2\fixTwo}
{\fixOne+\lambda^2\fixTwo}\Bigg|\by\right)\,\by.
$$
\label{res:threshSide}
\end{GHSresult}
A derivation of Result \ref{res:threshSide} is given in 
Section \ref{sec:derivResFive}.

For general Bayesian models containing (\ref{eq:bloodOrange}) or, equivalently,  
(\ref{eq:wetPubHol}) forms Result \ref{res:threshSide} suggests the following rule:
\begin{equation}
\mbox{decide that}\ \btheta=\bzero\quad\mbox{if and only if}\quad E(\gammaGHS|\by)< 
\smhalf\quad
\mbox{where}\quad 
\gammaGHS\equiv
\frac{\lambda^2\sigma_{\theta}^2}
{\sigeps^2+\lambda^2\sigma_{\theta}^2}.
\label{eq:BoomCrashOpera}
\end{equation}
To better understand the efficacy of (\ref{eq:BoomCrashOpera}), we ran a simulation
study similar to that in Section 4 of our recent article, 
He \myand Wand (2024), on generalized additive model selection.
The study involved the Bayesian generalized additive model given by equation (9) 
in He \myand Wand (2024) where $\dLin$ is the number of candidate predictors
that may have a zero or linear effect and $\dNon$ is the number of candidate
predictors that may have a zero, linear or non-linear effect.
The study involved both the set-up in He \myand Wand (2024), with the Laplace-Zero
and Grouped Lasso-Zero priors that are used in that article, and an alternative version 
with the likelihood taking the form
$$\by|\beta_0,\bbeta,
\bu_1,\ldots,\bu_{\dNon},\sigeps^2\sim 
N\left(\bone_n\beta_0+\bX\bbeta
+{\displaystyle\sum_{j=1}^{\dNon}}\bZ_j\bu_j,\sigeps^2\bI_n\right)
$$
where $\bbeta$ is a $(\dLin+\dNon)\times1$ vector of linear effects coefficients 
and, for each $1\le j\le\dNon$, $\bu_j$ is a $K_j\times1$ vector
of spline coefficients for the $j$th non-linear effect. 
Section 2 of He \myand Wand (2024) contains fuller details,
including the definition of the spline basis $\bZ_j$ matrices.

Let $\beta_j$ denote the $j$th entry of $\bbeta$.
Rather than imposing Laplace-Zero distributions on the $\beta_j$,
as conveyed by equation (6) of He \myand Wand (2024), we instead consider
the independent scalar Horseshoe specifications
\begin{equation}
\pDens(\beta_j|\sigma_{\beta})
=\pHSone(\beta_j/\sigma_{\beta})\big/\sigma_{\beta},\quad
1\le j\le\dLin+\dNon.
\label{eq;CayugaNights}
\end{equation}
Similarly, rather than imposing a Grouped Lasso-Zero distribution
on $\bu_j$, as conveyed by equation (7) of He \myand Wand (2024), we instead consider
the independent Grouped Horseshoe prior specifications
\begin{equation}
\pDens(\bu_j|\sigma_{u_j})
=\pHSKj(\bu_j/\sigma_{u_j})\Big/\sigma_{u_j}^{K_j},\quad
1\le j\le\dNon.
\label{eq:CollegetownDays}
\end{equation}
Note that (\ref{eq;CayugaNights}) has the auxiliary variable representation
$$\beta_j|\sigma_{\beta},\lambda_{\beta j}\simind
N(0,\sigma_{\beta}^2\lambda_{\beta j}),\quad 
\pDens(\lambda_{\beta j})=\frac{2I(\lambda_{\beta j}>0)}{\pi(1+\lambda_{\beta j}^2)},
\quad 1\le j\le\dLin+\dNon,$$
where $\simind$ denotes ``independently distributed as''.
Similarly, (\ref{eq:CollegetownDays}) has the auxiliary variable representation
$$\bu_j|\sigma_{uj},\lambda_{uj}\simind
N\big(\bzero,\sigma_{uj}^2\lambda_{uj}\bI_{K_j}\big),\quad 
\pDens(\lambda_{uj})=\frac{2I(\lambda_{uj}>0)}{\pi(1+\lambda_{uj}^2)},
\quad 1\le j\le\dNon.$$
From (\ref{eq:BoomCrashOpera}), the (Grouped) Horseshoe analogues of 
He \myand Wand (2024)'s $\gammaSUBbetaj$ and $\gammaSUBuj$ are
$$\gammabetajHS\equiv\frac{\lambda_{\beta_j}^2\sigma_{\beta}^2}
{\sigeps^2+\lambda_{\beta_j}^2\sigma_{\beta}^2}
\quad\mbox{and}\quad
\gammaujGHS\equiv\frac{\lambda_{uj}^2\sigma_{uj}^2}
{\sigeps^2+\lambda_{uj}^2\sigma_{uj}^2}.
$$
Therefore, for the $\dNon$ predictors that can have a zero, linear
or non-linear effect the classification rule that arises from 
(\ref{eq:BoomCrashOpera}) is
\begin{equation}
\begin{array}{l}
\mbox{the effect is zero if}\ 
\mbox{max}\big\{E(\gammabetajHS|\by),E(\gammaujGHS|\by)\big\}\le\smhalf,\\[1ex]
\mbox{the effect is linear if}\ E(\gammabetajHS|\by)>\smhalf\ \mbox{and}\
E(\gammaujGHS|\by)\le\smhalf,\\[1ex]
\mbox{otherwise the effect is non-linear.}
\end{array}
\label{eq:GHSclasRule}
\end{equation}
Note that this rule is analogous to the $\tau=\smhalf$ rule given in Section 3.5.2 of 
He \myand Wand (2024) for Laplace-Zero and Grouped Lasso-Zero priors,
with $\gammabetajHS$ and $\gammaujGHS$ instead of their
$\gammaSUBbetaj$ and $\gammaSUBuj$.

\begin{figure}[!t]
\centering
\includegraphics[width=0.95\textwidth]{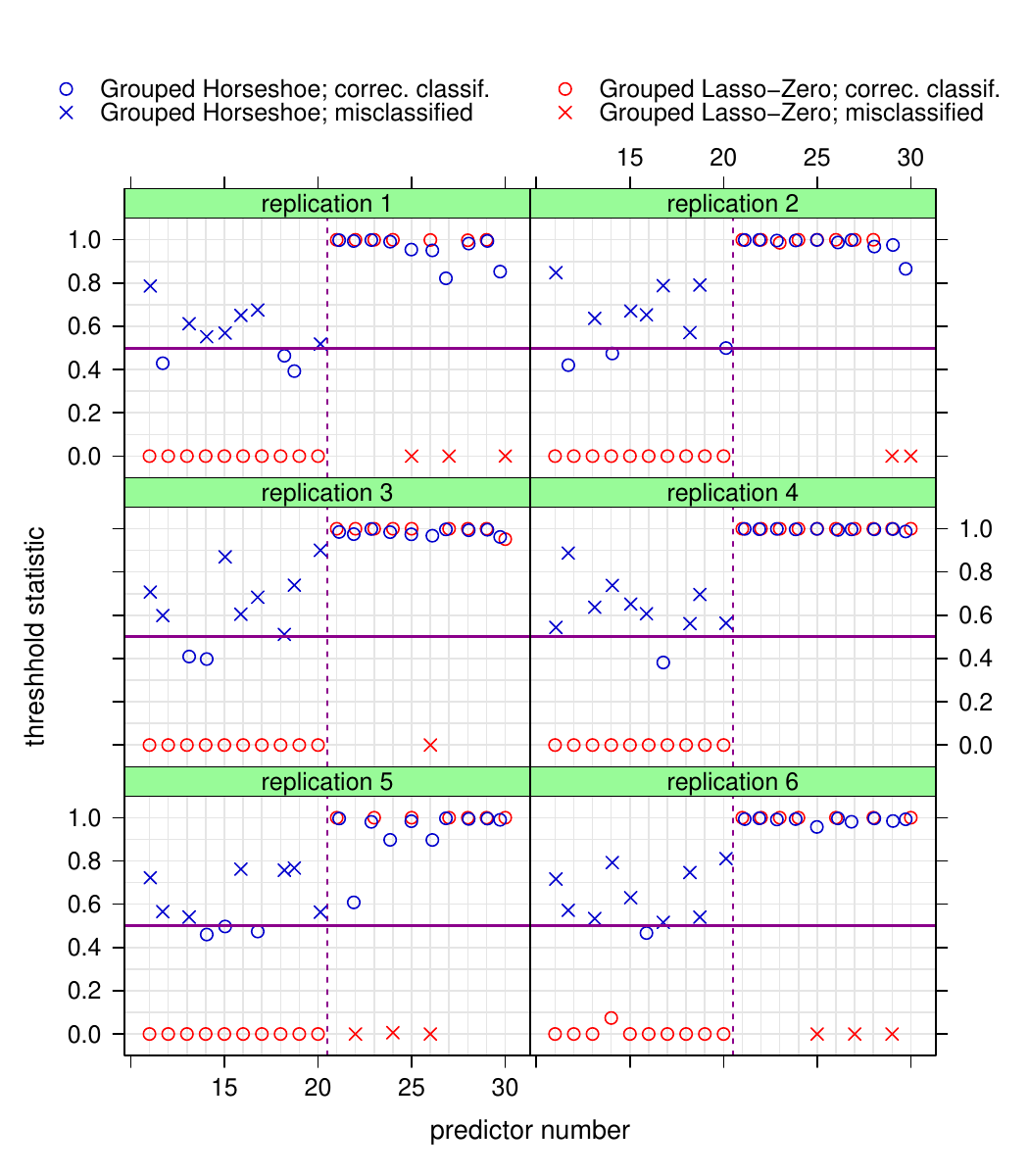}
\caption{
\textit{The results of linear effect versus non-linear effect classification
for six replications of the simulation study described in Section 4 of
He \myand Wand (2024) with $n=500$ and $\sigeps=2$.
The blue symbols correspond to the Markov chain Monte Carlo-approximate} 
$E(\gammaujGHS|\by)$ \textit{values for $21\le j\le30$,
which are the predictors that are simulated to have linear effects
(left of the vertical dashed line) or non-linear effects (right
of the vertical dashed line). A blue circle indicates correct classification
using (\ref{eq:GHSclasRule}), whilst a blue cross indicates misclassification.
The red symbols are similar, but for the $E(\gamma_{u_j}|\by)$ statistics
corresponding to strategy described in Section 3.5.2 of He \myand Wand (2024).
The classification border of $\smhalf$ is shown by the horizontal purple line.}}
\label{fig:GLZvsGHSclasComp1}
\end{figure}

Figure \ref{fig:GLZvsGHSclasComp1} shows the approximate,
based on Markov chain Monte Carlo sampling, values of 
$E(\gammaujGHS|\by)$ and $E(\gammaSUBuj|\by)$ for 
six replications of the simulation study described
in Section 4 of He \myand Wand (2024) 
with $n=500$ and $\sigeps=2$. The data are simulated
so that the effect of the $j$th predictor is
\begin{equation}
\begin{array}{ll}
\mbox{zero}       & \mbox{for}\quad j\in\{1,2,\ldots,10\},\\[0.1ex]
\mbox{linear}     & \mbox{for}\quad j\in\{11,12,\ldots,20\}\ \mbox{and}\\[0.1ex]
\mbox{non-linear} & \mbox{for}\quad j\in\{21,22,\ldots,30\}.
\end{array}
\label{eq:copperKettle}
\end{equation}
The horizontal axis corresponds to $j=11,12,\ldots,30$, which
is concerned with linear versus non-linear classification.
Correct classifications are shown as circles
and incorrect classifications are shown as crosses.
For these replications, use of the Grouped Lasso-Zero prior
results in a misclassification rate of $12/120=10\%$.
For the Grouped Horseshoe prior the misclassification
rate is $47/120=39.2\%$, which is about four times worse. 
We see from Figure \ref{fig:GLZvsGHSclasComp1}
that most of the Grouped Lasso-Zero threshold statistics 
are close to $1$ when the true effect is non-linear and close to 
$0$ when the true effect is linear. In contrast, most of the Grouped
Horseshoe threshold statistics are close to $1$ when the true effect is non-linear,
but scattered between $0.4$ and $0.8$ when the true effect is linear.
This last-mentioned behavior means that many predictors
that have a linear effect are misclassified as having
a non-linear effect when the Grouped Horseshoe prior is used.

\begin{figure}[!t]
\centering
\includegraphics[width=0.95\textwidth]{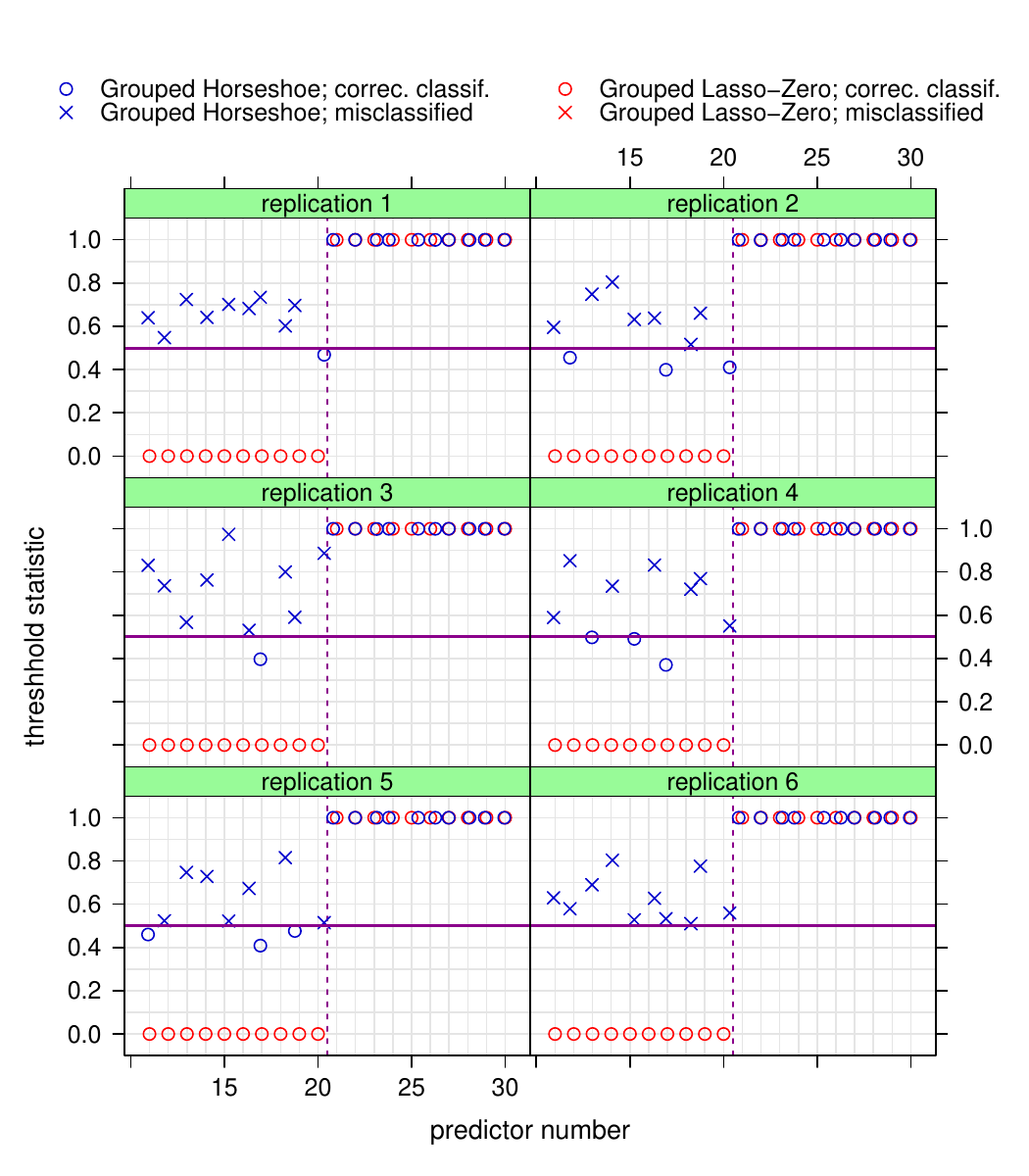}
\caption{
\textit{Similar to Figure \ref{fig:GLZvsGHSclasComp1}
but with a lower sample size, $n=2,000$, and a higher error standard
deviation, $\sigeps=0.25$.}}
\label{fig:GLZvsGHSclasComp2}
\end{figure}

Figure \ref{fig:GLZvsGHSclasComp2} differs from Figure \ref{fig:GLZvsGHSclasComp1}
in that the sample size is quadrupled to $n=2,000$ and the error standard deviation
is decreased to $\sigeps=0.25$. This should make linear versus
non-linear classification much easier and use of the Grouped Lasso-Zero
prior leads to perfect performance for these six replications.
However, for the Grouped Horseshoe prior the more favorable conditions
do not seem to help and the threshold statistics are still scattered
between $0.4$ and $0.8$ when the true effect is linear, 
leading to a $49/120=40.8\%$ misclassification rate.

\begin{figure}[!t]
\centering
\includegraphics[width=0.95\textwidth]{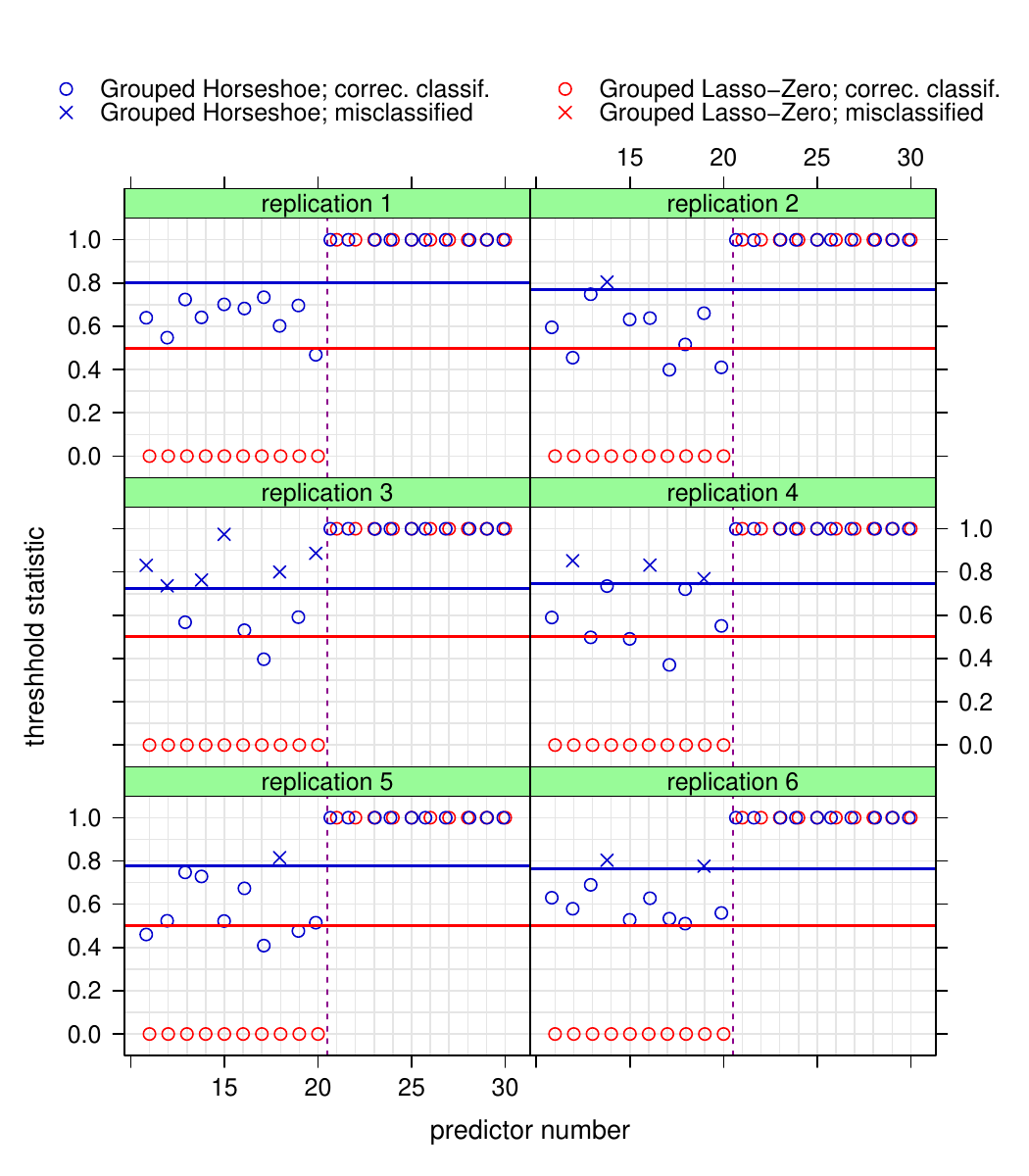}
\caption{
\textit{Similar to Figure \ref{fig:GLZvsGHSclasComp2}
but with the thresholding of the Grouped Horseshoe $E(\gammaujGHS|\by)$
statistics based on $k$-means clustering. For each replication,
the horizontal blue line shows the classification border arising
from $k$-means clustering. The horizontal red line at $\smhalf$ is 
the threshold for the Grouped Lasso-Zero statistics.}}
\label{fig:GLZvsGHSclasCompCLUS2}
\end{figure}

We experimented with a possible remedy to the poor performance of 
thresholding the $E(\gammaujGHS|\by)$ statistics at $\smhalf$.
This involved applying $k$-means clustering (e.g. MacQueen, 1967) 
to the $E(\gammaujGHS|\by)$ observations, 
with the number of clusters fixed at $2$. The function \texttt{kmeans()}
within the \textsf{R} computing environment (\textsf{R} Core Team, 2024)
was used to obtain the two clusters and corresponding classification rule.
As an example, for the analysis correponding to replication 1 of 
Figure \ref{fig:GLZvsGHSclasComp2}, the $k$-means threshold is $0.8031$. 
The results from use of this $k$-means alternative
to the Figure \ref{fig:GLZvsGHSclasComp2} analyses
are shown in Figure \ref{fig:GLZvsGHSclasCompCLUS2}. 
The misclassification rate drops to $13/120=10.83\%$.

This experimental $k$-means approach to thresholding for generalized additive model
selection has some promise, but relies on situations where there are
many candidate predictors of various effect types. If there are only 
$3$--$6$ candidate predictors, say, such that most of them have strongly 
non-linear effects then $k$-means threshold choice may not be viable.

Figures \ref{fig:GLZvsGHSclasComp1}--\ref{fig:GLZvsGHSclasCompCLUS2} are
based on only six replications. They also omit the
zero versus linear/non-linear classifications based on the 
$E(\gammabetajHS|\by)$ statistics and their Laplace-Zero counterparts,
which has similar results regarding Horseshoe versus Laplace-Zero priors.
To get a more complete picture, we ran an adaptation of the simulation
study described in Section 4 of He \myand Wand (2024). 
The data were generated in exactly the same manner as there, with
$30$ candidate predictors having ``true'' effects as described by (\ref{eq:copperKettle}),
the sample size ranging over $n\in\{500,1000,2000\}$ and the error standard
deviation ranging over $\sigeps\in\{0.25,0.5,1,2\}$.

\begin{figure}[h]
\centering
\includegraphics[width=0.98\textwidth]{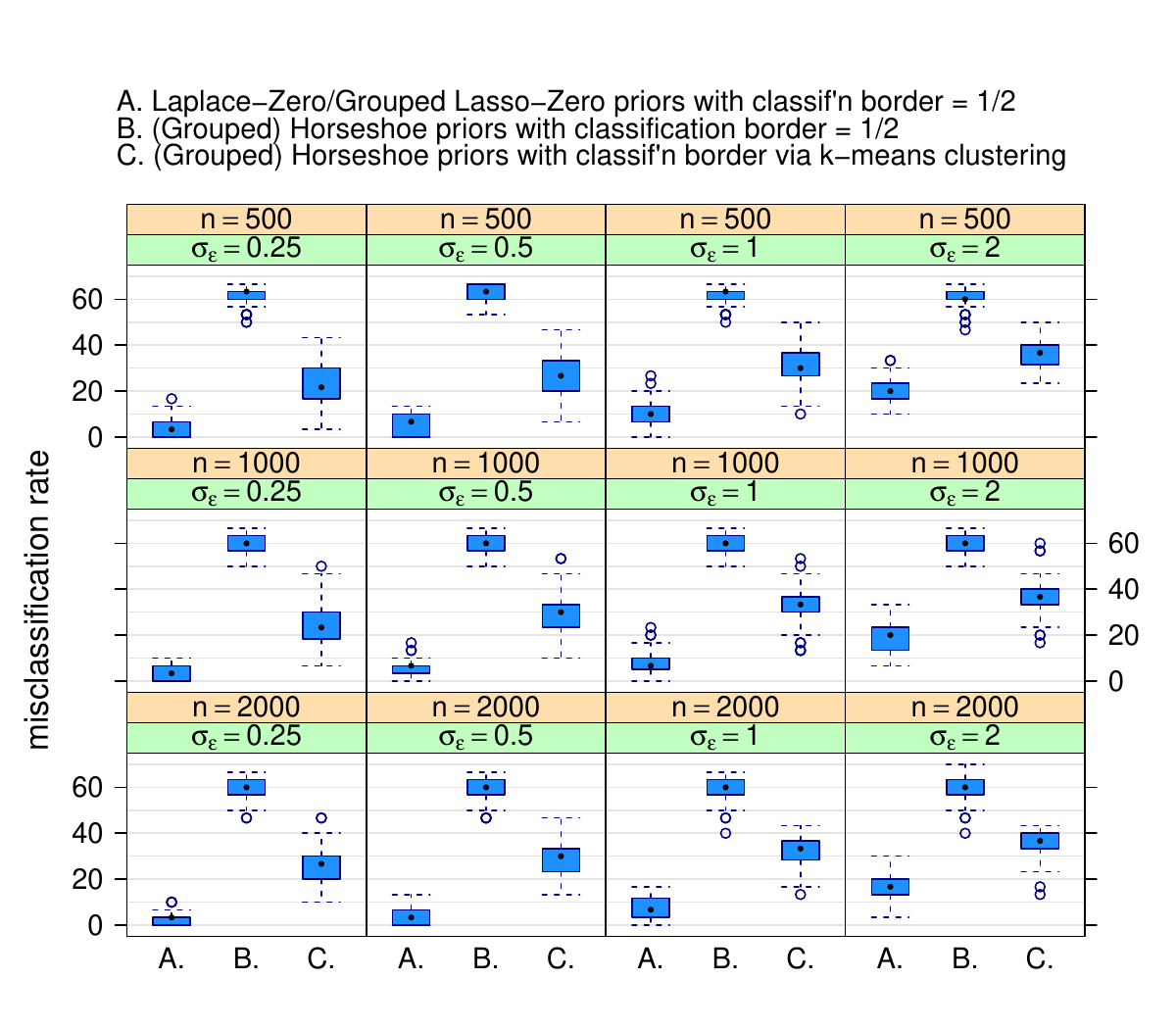}
\caption{
\textit{Side-by-side boxplots of the misclassification rates
for the comparative performance simulation study described in the text in the
case of the response variable being Gaussian. Each panel corresponds to a different 
combination of sample size and error standard deviation. Within each panel, the 
side-by-side boxplots compare the zero, linear, non-linear effect
misclassification rate across each of three methods:
A. Laplace-Zero/Grouped Lasso-Zero priors with classification border at $\smhalf$,
B. (Grouped) Horseshoe priors with classification border at $\smhalf$,
C. (Grouped) Horseshoe priors with classification border determined
via $k$-means clustering.}
}
\label{fig:GLZvsGHSsimRes}
\end{figure}

The side-by-side boxplots in Figure \ref{fig:GLZvsGHSsimRes} facilitate
comparison of 
\begin{itemize}
\item[A.] the Laplace-Zero/Grouped Lasso-Zero prior approach of He \myand Wand (2024)
with the classification border set to $\smhalf$,
\item[B.] use of the rule (\ref{eq:GHSclasRule}) involving 
(Grouped) Horseshoe priors, the $E(\gammabetajHS|\by)$ and
$E(\gammaujGHS|\by)$ threshold statistics and also with the 
classification border set to $\smhalf$,
\item[C.] the same as B., but with the classification border based on $k$-means clustering.
\end{itemize}
We see that A. is clearly superior to B. Also, C. offers a big improvement on B.,
but does not perform as well as A. In the He \myand Wand (2024) Bayesian generalized 
additive model selection setting use of (Grouped) Horseshoe priors does not compete
very well with use of Laplace-Zero/Grouped Lasso-Zero priors.

\section{Conclusions}\label{sec:conclusions}

In this article we have conducted a thorough investigation into
the statistical properties of the Grouped Horseshoe distribution.
We have shown that most of the properties possessed by the 
univariate Horseshoe distribution extend to the grouped situation.
Our investigation was motivated by our interest in Bayesian generalized
additive model selection, as described in our recent He \myand Wand (2024)
article. The numerical studies in Section \ref{sec:thresh}, concerned with
using Result \ref{res:threshSide} to carry out generalized
additive model selection with (Grouped) Horseshoe priors, 
reveal some performance concerns compared with a spike-and-slab benchmark.
Perhaps this research can lead to the development of better 
selection rules based on (Grouped) Horseshoe priors.

\section*{Acknowledgements}

We are grateful to Andrew Barron, Anindya Bhadra and Marty Wells 
for their contributions to this research. This research was partially 
supported by Australian Research Council grant DP230101179. 
The second author is grateful for hospitality from the 
Department of Statistics and Data Science, Cornell University, U.S.A., during 
part of this research.

\section*{References}

\bib
Bai, R. \myand Ghosh, M. (2018).
High dimensional multivariate posterior consistency under global-local
shrinkage priors. \textit{Journal of Multivariate Analysis}, 
\textbf{167}, 157--170.

\bib
Carvalho, C.M., Polson, N.G. \myand Scott, J.G. (2010).
The horseshoe estimator for sparse signals.
\textit{Biometrika}, \textbf{97}, 465--480.

\bib
He, V.X. and Wand, M.P. (2024).
Bayesian generalized additive model selection including
a fast variational option.
\textit{Advances in Statistical Analysis}, in press.

\bib
MacQueen, J.B. (1967). Some methods for classification and analysis
of multivariate observations. In \textit{Proceedings of 
the Fifth Berkeley Symposium on Mathematical Statistics 
and Probability, Volume 1}, pp. 281--297, Berkeley, California:
\textit{University of California Press}.

\bib
\textsf{R} Core Team (2024). \textsf{R}: A language and environment for statistical
computing. \textsf{R} Foundation for Statistical Computing, Vienna, Austria.
\texttt{https://www.R-project.org/}

\bib
Scheipl, F., Fahrmeir, L. \myand Kneib, T. (2012). 
Spike-and-slab priors for function selection in structured
additive regression models.\ {\it Journal of the American
Statistical Association}, {\bf 107}, 1518--1532.

\bib
Xu, Z., Schmidt, D.F., Enes, M., Qian, G. \myand Hopper, J.L. (2016).
Bayesian grouped horseshoe regression with application
to additive models.
In B.H. Kang \myand Q. Bai (eds.), \textit{AI 2016: Advances
in Artificial Intelligence}, pp. 229--240,
Cham, Switzerland: Springer.

\null
\vfill\eject
%
%
\renewcommand{\theequation}{S.\arabic{equation}}
\renewcommand{\thesection}{S.\arabic{section}}
\renewcommand{\thetable}{S.\arabic{table}}
\setcounter{equation}{0}
\setcounter{table}{0}
\setcounter{section}{0}
\setcounter{page}{1}
\setcounter{footnote}{0}

\centerline{\Large Supplement for:}
\vskip5mm
\centerline{\Large\bf The Grouped Horseshoe Distribution and Its Statistical Properties}
\vskip5mm
\centerline{\normalsize\sc Virginia X. He and Matt P. Wand}
\vskip5mm
\centerline{\textit{University of Technology Sydney}}
\vskip5mm{\null}

\section{Derivation of Result \ref{res:explicDens}}\label{sec:derivResOne}

From (\ref{eq:GilbertGrape}) we have 
$$\pHSd(\bx)=\int_0^{\infty}(2\pi\lambda^2)^{-d/2}
\exp\left(-\frac{\Vert\bx\Vert^2}{2\lambda^2}\right)\,
\frac{2}{\pi(1+\lambda^2)}\,d\lambda.
$$
The change of variable $t=1/\lambda^2$ then leads to 
{\setlength\arraycolsep{1pt}
\begin{eqnarray*}
\pHSd(\bx)&=&(2\pi)^{-d/2}\pi^{-1}\int_0^{\infty}
\frac{t^{(d+1)/2-1}\exp(-t\Vert\bx\Vert^2/2)}{1+t}\,dt\\[1ex]
&=&(2\pi)^{-d/2}\pi^{-1}\Gamma\big(\smhalf(d+1)\big)\big(\Vert\bx\Vert^2/2\big)^{(1-d)/2}
\exp\big(\Vert\bx\Vert^2/2\big)\\[1ex]
&&\quad\times\frac{\big(\Vert\bx\Vert^2/2\big)^{(d+1)/2-1}}
{\Gamma\big(\smhalf(d+1)\big)}\exp\big(-\Vert\bx\Vert^2/2\big)
\int_0^{\infty}
\frac{t^{(d+1)/2-1}\exp(-t\Vert\bx\Vert^2/2)}{1+t}\,dt\\[1ex]
&=&\frac{\Gamma\big(\smhalf(d+1)\big)}{\sqrt{2\pi^{d+2}}}\,
\exp\big(\Vert\bx\Vert^2/2\big)
E_{(d+1)/2}\big(\Vert\bx\Vert^2/2\big)\Big/\Vert\bx\Vert^{d-1}
\end{eqnarray*}
}
with the last step following from 8.19.4 of Olver (2023).

\section{Derivation of Result \ref{res:poleExist}}\label{sec:derivResTwo}

We break up the derivation into the cases:
$$d=1\quad\mbox{and}\quad d\ge2.$$

\vskip2mm
\leftline{$\underline{\mbox{\textit{The $d=1$ Case}}}$}
\vskip2mm

As stated in Theorem 1 of Carvalho \textit{et al.} (2010),
$$\pHSone(x)>\frac{K_1}{2}\log\left(1+\frac{4}{x^2}\right)\quad\mbox{where}\quad 
K_1\equiv \frac{1}{\sqrt{2\pi^3}}.$$
Then
$$\lim_{x\to 0}\pHSone(x)>
\frac{K_1}{2}\log\left(1+4\lim_{x\to0}\left(\frac{1}{x^2}\right)\right)
=\infty.$$
Hence 
$$\lim_{x\to 0}\pHSone(x)=\infty$$
and Result \ref{res:poleExist} holds for $d=1$.

\vfill\eject
\leftline{$\underline{\mbox{\textit{The $d\ge2$ Case}}}$}
\vskip2mm

From Result \ref{res:explicDens}
$$\pHSd(\bx)=K_d\,\exp\big(\Vert\bx\Vert^2/2\big)
E_{(d+1)/2}\big(\Vert\bx\Vert^2/2\big)\Big/\Vert\bx\Vert^{d-1}
\quad\mbox{where}\quad K_d\equiv
\frac{\Gamma\big(\smhalf(d+1)\big)}{\sqrt{2\pi^{d+2}}}.$$
Then 
$$\lim_{\bx\to\bzero}\pHSd(\bx)=K_d\left\{\lim_{\bx\to\bzero} 
\exp\big(\Vert\bx\Vert^2/2\big) \right\}
\left\{\lim_{\bx\to\bzero} 
E_{(d+1)/2}\big(\Vert\bx\Vert^2/2\big)\right\}
\left\{\lim_{\bx\to\bzero}\Vert\bx\Vert^{1-d}\right\}.
$$
Clearly, 
$$\lim_{\bx\to\bzero}\exp\big(\Vert\bx\Vert^2/2\big)=1.$$
Also, from 8.19.6 of Olver (2023),
$$E_{\nu}(0)=\frac{1}{\nu-1},\quad\mbox{for all}\quad \nu>1,$$
which leads to
$$\lim_{\bx\to\bzero} 
E_{(d+1)/2}\big(\Vert\bx\Vert^2/2\big)
=\frac{1}{\smhalf(d+1)-1}=\frac{2}{d-1}\in(0,2]\quad
\mbox{for all}\quad d\ge 2.
$$
Lastly, 
$$\lim_{\bx\to\bzero}\Vert\bx\Vert^{1-d}=\infty\quad\mbox{for all}\quad d\ge2.$$
Hence 
$$\lim_{\bx\to\bzero}\pHSd(\bx)=\infty\quad\mbox{for all}\quad d\ge2.$$

\section{Derivation of Result \ref{res:robustRes}}\label{sec:derivResThree}

Result \ref{res:robustRes} entails various properties of the special function
known as the bivariate confluent hypergeometric function. We commence with 
its definition and some key results. We then show how these results lead to 
the Result \ref{res:robustRes} statements.

\subsection{The Bivariate Confluent Hypergeometric Function}\label{sec:BCHF}

The \emph{bivariate confluent hypergeometric function} 
$$\Phi_1(\alpha,\beta,\gamma,x,y)\quad\mbox{for}\quad\alpha,\beta,\gamma,x,y\in{\mathbb C}$$
is defined via a pair of partial differential equations in Section 9.262
of Gradshteyn \myand Ryzhik (1994). As stated in Section 9.261 of 
Gradshteyn \myand Ryzhik (1994) the following series representation applies 
when $|x|<1$:
\begin{equation}
\Phi_1(\alpha,\beta,\gamma,x,y)=\sum_{m=0}^{\infty}\sum_{n=0}^{\infty}
\frac{(\alpha)_{m+n}(\beta)_m x^m y^n}{(\gamma)_{m+n}m!n!}\quad\mbox{where}
\quad(a)_k\equiv\Gamma(a+k)/\Gamma(a).
\label{eq:BCHFseries}
\end{equation}

Result 3.385 of Gradshteyn \myand Ryzhik (1994) states that
\begin{equation}
\int_0^1 x^{\nu-1}(1-x)^{\lambda-1}(1-\beta x)^{-\rho}e^{-\mu x}\,dx
=\frac{\Gamma(\nu)\Gamma(\lambda)}{\Gamma(\nu+\lambda)}\,
\Phi_1(\nu,\rho,\nu+\lambda,\beta,-\mu)
\label{eq:GRinteg}
\end{equation}
for complex numbers $\lambda$, $\nu$, $\rho$, $\beta$ and $\mu$ ranging
over various subsets of the complex plane. If these parameters are constrained
to be real then the restrictions reduce to 
$$\lambda,\nu>0\quad\mbox{and}\quad\beta,\rho,\mu\in\real.$$

Arguments in Appendix B of Gordy (1998) imply that for $0\le x<1$ and $0<\alpha<\gamma$ 
we have the following series representation of $\Phi_1$ in terms of the 
univariate confluent hypergeometric function $\oneFone$:
\begin{equation}
\Phi_1(\alpha,\beta,\gamma,x,y)
=\exp(y)\sum_{n=0}^{\infty}
\frac{(\alpha)_n(\beta)_n x^n}{(\gamma)_n\,n!}\,
\oneFone(\gamma-\alpha,\gamma+n,-y).
\label{eq:PhiOneOneFone}
\end{equation}
Note, however, that there in an error in equation (6) of Gordy (1998). 
It is due to the $(\beta)_m$ of (\ref{eq:BCHFseries}) being incorrectly 
replaced by $(\beta)_n$.
This error leads to (T1) -- (T4) of Gordy (1998)
containing an incorrect variant of (\ref{eq:PhiOneOneFone}) with
respect to the $\Phi_1$, $\oneFone$ and $\twoFone$ functions 
as defined in Gradshteyn \myand Ryzhik (1994).

\subsection{Marginal Density Function Simplification}

The marginal density function of $\by$ according to model (\ref{eq:robMod}) is
$$\pDens(\by)=\int_0^{\infty}\left\{\int_{\real^d}\pDens(\by|\btheta)
\pDens(\btheta|\lambda)d\btheta\right\}\pDens(\lambda)\,d\lambda\quad\mbox{where}\quad
\pDens(\lambda)=\frac{2I(\lambda>0)}{\pi(\lambda^2+1)}.
$$
Define
$$\kappaREPL\equiv 1 + 1/(\lambda^2\tau^2).$$
Then standard algebraic arguments lead to 
{\setlength\arraycolsep{1pt}
\begin{eqnarray*}
\pDens(\by|\btheta) \pDens(\btheta|\lambda)
&=&(2\pi)^{-d/2}(\lambda\tau)^{-d}\exp[\smhalf\{(1/\kappaREPL)-1\}\Vert\by\Vert^2]|(1/\kappaREPL)\bI_d|^{1/2}\\
&&\ \times
(2\pi)^{-d/2}|(1/\kappaREPL)\bI_d|^{-1/2}\exp\big[-\smhalf\{\btheta-(1/\kappaREPL)\by\}^T\{(1/\kappaREPL)\bI_d\}^{-1}
\{\btheta-(1/\kappaREPL)\by\}\big].
\end{eqnarray*}
}
Noting that
$$(2\pi)^{-d/2}|(1/\kappaREPL)\bI_d|^{-1/2}\exp\big[-\smhalf\{\btheta-(1/\kappaREPL)\by\}^T\{(1/\kappaREPL)\bI_d\}^{-1}
\{\btheta-(1/\kappaREPL)\by\}\big]$$
is the $N\Big((1/\kappaREPL)\by,(1/\kappaREPL)\bI_d\Big)$ density function in $\btheta$
and $|(1/\kappaREPL)\bI_d|=\kappaREPL^{-d}$, we then have
{\setlength\arraycolsep{1pt}
\begin{eqnarray*}
\int_{\real^d}\pDens(\by|\btheta)\pDens(\btheta|\lambda)d\btheta
&=&(2\pi)^{-d/2}(\lambda\tau)^{-d}\exp[\smhalf\{(1/\kappaREPL)-1\}\Vert\by\Vert^2]
\kappaREPL^{-d/2}\\[1ex]
&=&(2\pi)^{-d/2}\exp\left\{-\frac{(\Vert\by\Vert^2/2)}{1+\lambda^2\tau^2}\right\}
\frac{1}{(1+\lambda^2\tau^2)^{d/2}}.
\end{eqnarray*}
}
The marginal density function of $\by$ is then
\begin{equation}
\pDens(\by)=(2^{d-2}\pi^{d+2})^{-1/2}\int_0^{\infty}\
\exp\left\{-\frac{(\Vert\by\Vert^2/2)}{1+\lambda^2\tau^2}\right\}
\frac{1}{(1+\lambda^2\tau^2)^{d/2}(\lambda^2+1)}\,d\lambda.
\label{eq:Laramie}
\end{equation}

\subsection{Score Function Simplification}

The score function is the following $d\times 1$ derivative vector:
$$\nabla_{\by}\{\log \pDens(\by)\}.$$
Next note that 
{\setlength\arraycolsep{1pt}
\begin{eqnarray*}
d_{\by}\log \pDens(\by)&=&\frac{1}{\pDens(\by)}d_{\by}\,\pDens(\by)\\[1ex]
&=&\displaystyle{\frac{\displaystyle{\int_0^{\infty}\
d_{\by}\exp\left\{-\frac{(\Vert\by\Vert^2/2)}{1+\lambda^2\tau^2}\right\}
\frac{1}{(1+\lambda^2\tau^2)^{d/2}(\lambda^2+1)}\,d\lambda}}
{\displaystyle{\int_0^{\infty}\
\exp\left\{-\frac{(\Vert\by\Vert^2/2)}{1+\lambda^2\tau^2}\right\}
\frac{1}{(1+\lambda^2\tau^2)^{d/2}(\lambda^2+1)}\,d\lambda}}}\\[1ex]
&=&-\by\left[\displaystyle{\frac{
\displaystyle{\int_0^{\infty}
\exp\left\{-\frac{(\Vert\by\Vert^2/2)}{1+\lambda^2\tau^2}\right\}
\frac{1}{(1+\lambda^2\tau^2)^{(d+2)/2}(\lambda^2+1)}\,d\lambda}}
{\displaystyle{\int_0^{\infty}\
\exp\left\{-\frac{(\Vert\by\Vert^2/2)}{1+\lambda^2\tau^2}\right\}
\frac{1}{(1+\lambda^2\tau^2)^{d/2}(\lambda^2+1)}\,d\lambda}}}\right]\,d\by.
\end{eqnarray*}
}
Hence 
\begin{equation}
\nabla_{\by}\{\log \pDens(\by)\}=
-\by\left[\displaystyle{\frac{
\displaystyle{
\int_0^{\infty}
\exp\left\{-\frac{(\Vert\by\Vert^2/2)}{1+\lambda^2\tau^2}\right\}
\frac{1}{(1+\lambda^2\tau^2)^{(d+2)/2}(\lambda^2+1)}\,d\lambda}}
{\displaystyle{\int_0^{\infty}\
\exp\left\{-\frac{(\Vert\by\Vert^2/2)}{1+\lambda^2\tau^2}\right\}
\frac{1}{(1+\lambda^2\tau^2)^{d/2}(\lambda^2+1)}\,d\lambda}}}\right].
\label{eq:CubsCap}
\end{equation}

\subsection{Bivariate Confluent Hypergeometric Function Representations}\label{sec:BCHFreps}

In this subsection we derive expressions for $\pDens(\by)$ and 
$\nabla_{\by}\{\log \pDens(\by)\}$ in terms of the
bivariate confluent hypergeometric function $\Phi_1$ 
as defined in Section \ref{sec:BCHF}. 

The integral in (\ref{eq:Laramie}) is 
$$\Csc\big(\smhalf\Vert\by\Vert^2,\tau^2\big)$$
where
\begin{equation}
\Csc(a,b)\equiv \int_0^{\infty}\
\exp\left(-\frac{a}{1+\lambda^2b}\right)
\frac{1}{(1+\lambda^2b)^{d/2}(\lambda^2+1)}\,d\lambda.
\label{eq:TalkingMule}
\end{equation}
The change of variable 
$$x=\lambda^2b/(1+\lambda^2b)$$
in the (\ref{eq:TalkingMule}) integral leads to 
$$
\Csc(a,b)=\frac{\exp(-a)}{2\sqrt{b}}\int_0^1 x^{\nu-1} (1-x)^{\lambda-1}
(1-\beta x)^{-\rho}e^{-\mu\,x}\,dx$$
where 
$$\nu=\smhalf,\quad\lambda=\smhalf(d+1),\quad\beta=1-b^{-1},\quad\rho=1
\quad\mbox{and}\quad\mu=-a.
$$
Application of (\ref{eq:GRinteg}) provides the bivariate confluent hypergeometric
form
$$\Csc(a,b)=\frac{\exp(-a)\sqrt{\pi}\Gamma\big(\smhalf d+\smhalf\big)}
{d\sqrt{b}\,\Gamma\big(\smhalf d\big)}
\Phi_1\Big(\smhalf,1,\smhalf(d+2),1-b^{-1},a\Big).
$$
Plugging this into (\ref{eq:Laramie}), with $a=\smhalf\Vert\by\Vert^2$ and 
$b=\tau^2$, we obtain
\begin{equation}
\pDens(\by)=\big(2^{d-2}\pi^{d+1}\big)^{-1/2}
\frac{\exp\big(-\smhalf\Vert\by\Vert^2\big)\Gamma\big(\smhalf d+\smhalf\big)}
{\tau\,d\,\Gamma\big(\smhalf d\big)}\,
\Phi_1\left(\smhalf,1,\smhalf(d+2),1-\tau^{-2},\smhalf\Vert\by\Vert^2\right).
\label{eq:pyBCHF}
\end{equation}
For the $d=1$ special case, (\ref{eq:pyBCHF}) reduces to an expression similar, but not
identical, to that provided by equation (A1) of Carvalho \textit{et al.} (2010). 
The main difference is an interchange in the fourth and fifth arguments of the $\Phi_1$ function.
This discrepancy is attributable to an error in Gordy (1998), which we  
described in Section \ref{sec:BCHF}.

Next we seek an analogous expression  for the score function.
It follows from (\ref{eq:CubsCap}) that
$$\nabla_{\by}\log\{\pDens(\by)\}=\,-\by\frac{\Dsc\big(\smhalf\Vert\by\Vert^2,\tau^2\big)}
{\Csc\big(\smhalf\Vert\by\Vert^2,\tau^2\big)}
$$
where
$$\Dsc(a,b)\equiv \int_0^{\infty}\
\exp\left(-\frac{a}{1+\lambda^2b}\right)
\frac{1}{(1+\lambda^2b)^{(d/2)+1}(\lambda^2+1)}\,d\lambda.
$$
Calculations similar to those given in the previous section lead to 
$$\Dsc(a,b)=\frac{\exp(-a)}{2\sqrt{b}}\int_0^1 x^{\nu-1} (1-x)^{\lambda-1}
(1-\beta x)^{-\rho}e^{-\mu\,x}\,dx.
$$
where
$$\nu=\smhalf,\quad\lambda=\smhalf(d+3),\quad\beta=1-b^{-1},\quad\rho=1
\quad\mbox{and}\quad\mu=-a.
$$
From (\ref{eq:GRinteg}), 
$$\Dsc(a,b)=\frac{\exp(-a)}{2\sqrt{b}}\frac{\sqrt{\pi}\Gamma\big(\smhalf d+\threehalves\big)}
{\Gamma\big(\smhalf\,d+2\big)}\Phi_1\Big(\smhalf,1,\smhalf(d+4),1-b^{-1},a\Big).
$$
Next note that 
$$\Gamma\big(\smhalf d+\threehalves\big)=\Gamma\big(\smhalf d+\smhalf+1\big)=\smhalf(d+1)
\Gamma\big(\smhalf d+\smhalf\big)$$
and
$$\Gamma\big(\smhalf d+2\big)=\Gamma\big(\smhalf d+1+1\big)=(\smhalf d+1)
\Gamma\big(\smhalf d+1\big).$$
This leads to 
$$\frac{\Dsc(a,b)}{\Csc(a,b)}=\frac{(d+1)\Phi_1\Big(\smhalf,1,\smhalf(d+4),1-b^{-1},a\Big)}
{(d+2)\Phi_1\Big(\smhalf,1,\smhalf(d+2),1-b^{-1},a\Big)}.
$$
Hence
\begin{equation}
\nabla_{\by}\log\{\pDens(\by)\}=\,-
\frac{\by(d+1)\Phi_1\left(\smhalf,1,\smhalf(d+4),1-\tau^{-2},\smhalf\Vert\by\Vert^2\right)}
{(d+2)\Phi_1\Big(\smhalf,1,\smhalf(d+2),1-\tau^{-2},\smhalf\Vert\by\Vert^2\Big)}.
\label{eq:suxEggs}
\end{equation}
When $d=1$ this result matches equation (A2) of Carvalho \textit{et al.} (2010) 
except for interchanges in the fourth and fifth arguments of the 
$\Phi_1$ function. An error in Gordy (1998), which is described in 
Section \ref{sec:BCHF}, provides an explanation for this discrepancy.

\subsection{Large $\Vert\by\Vert$ Approximation of the Score Function}

In keeping with Result \ref{res:robustRes} being concerned with the limiting
behaviour of the score function as $\Vert\by\Vert\to\infty$, throughout
this subsection we assume that $\Vert\by\Vert\gg 1$. The following cases
are treated separately (in order of complexity):
$$\tau=1,\quad\tau>1\quad\mbox{and}\quad0<\tau<1.$$
In each case versions of the following result, from 
e.g. Section 13.1.5 of Abramowitz \myand Stegun (1968),
concerning the right-tail asymptotic behaviour of the 
univariate confluent hypergeometric function:
\begin{equation}
\oneFone(a,b,x)=\left\{
\begin{array}{lcr}
\displaystyle{\frac{\Gamma(b)}{\Gamma(a)}}
\,e^x x^{a-b}\big\{1+O(|x|^{-1})\big\} & \mbox{for}\ x>0\ \mbox{and}\ |x|\gg1,\\[2ex]
\displaystyle{\frac{\Gamma(b)}{\Gamma(b-a)}}
(-x)^{-a}\big\{1+O(|x|^{-1})\big\} & \mbox{for}\ x<0\ \mbox{and}\ |x|\gg1.\\
\end{array}
\right.
\label{eq:StegunRes}
\end{equation}

\vfill\eject
\leftline{$\underline{\mbox{\textit{The $\tau=1$ Case}}}$}
\vskip2mm

If $\tau=1$ it follows from Section 3.383 of Gradshteyn \myand Ryzhik (1994)
that
\begin{equation}
\nabla_{\by}\log\{\pDens(\by)\}=\,-
\frac{\by(d+1)\oneFone\left(\smhalf,\smhalf(d+4),\smhalf\Vert\by\Vert^2\right)}
{(d+2)\oneFone\Big(\smhalf,\smhalf(d+2),\smhalf\Vert\by\Vert^2\Big)}.
\label{eq:jjjSuz}
\end{equation}
From (\ref{eq:StegunRes}),
$$\oneFone\big(\smhalf,\smhalf(d+4),\smhalf\Vert\by\Vert^2\big)=
\frac{\Gamma\big(\smhalf(d+4)\big)}{\Gamma(\smhalf)}
\exp\big(\smhalf\Vert\by\Vert^2\big)\big(\smhalf\Vert\by\Vert^2\big)^{-\smhalf(d+3)}
\big\{1+O(\Vert\by\Vert^{-2})\big\}
$$
and
$$\oneFone\big(\smhalf,\smhalf(d+2),\smhalf\Vert\by\Vert^2\big)=
\frac{\Gamma\big(\smhalf(d+2)\big)}{\Gamma(\smhalf)}
\exp\big(\smhalf\Vert\by\Vert^2\big)\big(\smhalf\Vert\by\Vert^2\big)^{-\smhalf(d+1)}
\big\{1+O(\Vert\by\Vert^{-2})\big\}.
$$
Therefore, 
{\setlength\arraycolsep{1pt}
\begin{eqnarray*}
\frac{\oneFone\big(\smhalf,\smhalf(d+4),\smhalf\Vert\by\Vert^2\big)}
{\oneFone\big(\smhalf,\smhalf(d+2),\smhalf\Vert\by\Vert^2\big)}
&=&
\frac{2\big\{\Gamma\big(\smhalf(d+4)\big)/\Gamma\big(\smhalf(d+2)\big)\big\}
\Vert\by\Vert^{-2}\{1+O(\Vert\by\Vert^{-2})\}}
{1+O(\Vert\by\Vert^{-2})}\\[1ex]
&=&
\frac{(d+2)
\Vert\by\Vert^{-2}\{1+O(\Vert\by\Vert^{-2})\}}
{1+O(\Vert\by\Vert^{-2})}.
\end{eqnarray*}
}
We then have
\begin{equation}
\nabla_{\by}\log\{\pDens(\by)\}=\,-
\frac{(d+1)
\by\Vert\by\Vert^{-2}\{1+O(\Vert\by\Vert^{-2})\}}
{1+O(\Vert\by\Vert^{-2})}
\quad\mbox{for}\quad \Vert\by\Vert\gg1.
\label{eq:tauEQ1}
\end{equation}

\vskip2mm
\leftline{$\underline{\mbox{\textit{The $\tau> 1$ Case}}}$}
\vskip2mm

If $\tau>1$ then $0< 1-\tau^{-2}<1$ and we can use (\ref{eq:PhiOneOneFone}) to obtain
\begin{eqnarray*}
&&\Phi_1\left(\smhalf,1,\smhalf(d+4),1-\tau^{-2},\smhalf\Vert\by\Vert^2\right)\\
&&\qquad\qquad
=\exp\big(\smhalf\Vert\by\Vert^2\big)\sum_{n=0}^{\infty}
\frac{(\smhalf)_n(1-\tau^{-2})^n}{\big(\smhalf(d+4)\big)_n}\,
\oneFone\Big(\smhalf(d+3),\smhalf(d+4)+n,-\smhalf\Vert\by\Vert^2\Big).
\end{eqnarray*}
Next, from (\ref{eq:StegunRes}),
$$\oneFone\Big(\smhalf(d+3),\smhalf(d+4)+n,-\smhalf\Vert\by\Vert^2\Big)
=\frac{\Gamma\big(\smhalf(d+4)+n\big)}{\Gamma(n+\smhalf)}\{\smhalf\Vert\by\Vert^2\}^{-(d+3)/2}
\big\{1+O(\Vert\by\Vert^{-2})\big\}
$$
which leads to
$$\Phi_1\left(\smhalf,1,\smhalf(d+4),1-\tau^{-2},\smhalf\Vert\by\Vert^2\right)
=\frac{\tau^22^{(d+3)/2}\Gamma(\smhalf d+2)}{\sqrt{\pi}}\,
\exp\big(\smhalf\Vert\by\Vert^2\big)\Vert\by\Vert^{-(d+3)}
\big\{1+O(\Vert\by\Vert^{-2})\big\}.
$$
Similarly,
$$\Phi_1\left(\smhalf,1,\smhalf(d+2),1-\tau^{-2},\smhalf\Vert\by\Vert^2\right)
=\frac{\tau^22^{(d+1)/2}\Gamma(\smhalf d+1)}{\sqrt{\pi}}\,
\exp\big(\smhalf\Vert\by\Vert^2\big)\Vert\by\Vert^{-(d+1)}
\big\{1+O(\Vert\by\Vert^{-2})\big\}.
$$
Substitution into (\ref{eq:suxEggs}) then leads to 
\begin{equation}
\nabla_{\by}\log\{\pDens(\by)\}=\frac{-(d+1)\by\Vert\by\Vert^{-2}
\big\{1+O\big(\Vert\by\Vert^{-2}\big)\big\}}{1+O\big(\Vert\by\Vert^{-2}\big)}
\quad\mbox{for}\quad \Vert\by\Vert\gg1.
\label{eq:tauGT1}
\end{equation}
%

\vfill\eject
\leftline{$\underline{\mbox{\textit{The $0<\tau<1$ Case}}}$}
\vskip2mm

Recall from results in Section \ref{sec:BCHFreps} that
\begin{equation}
\nabla_{\by}\log\{\pDens(\by)\}=\,-\by\frac{\Dsc\big(\smhalf\Vert\by\Vert^2,\tau^2\big)}
{\Csc\big(\smhalf\Vert\by\Vert^2,\tau^2\big)}
\label{eq:goodEgg}
\end{equation}
where
$$
\Csc(a,b)=\frac{\exp(-a)}{2\sqrt{b}}\int_0^1 x^{-\smhalf} (1-x)^{\smhalf(d-1)}
\{1- (1-b^{-1}) x\}^{-1}e^{a\,x}\,dx$$
and
$$\Dsc(a,b)=\frac{\exp(-a)}{2\sqrt{b}}\int_0^1 x^{-\smhalf} (1-x)^{\smhalf(d+1)}
\{1-(1-b^{-1})x\}^{-1}e^{a\,x}\,dx.
$$
Noting that
$$1- (1-b^{-1})x=b^{-1}\{1-(1-b)(1-x)\}$$
we have
$$
\Csc(a,b)=\smhalf\exp(-a)\sqrt{b}\int_0^1 x^{-\smhalf} (1-x)^{\smhalf(d-1)}
\{1-(1-b)(1-x)\}^{-1}e^{a\,x}\,dx.
$$
The change of variables $u=1-x$ leads to
$$
\Csc(a,b)=\smhalf\sqrt{b}\int_0^1 u^{\nu-1} (1-u)^{\lambda-1} 
(1-\beta u)^{-\rho}e^{-\mu u}\,du
$$
where
$$\nu=\smhalf(d+1),\quad\lambda=\smhalf,\quad\beta=1-b,\quad\rho=1
\quad\mbox{and}\quad\mu= a.
$$
Hence, from (\ref{eq:GRinteg}),
$$\Csc(a,b)=\frac{\sqrt{b\pi}\,\Gamma\big(\smhalf(d+1)\big)}
{2\Gamma\big(\smhalf(d+2)\big)}\Phi_1\Big(\smhalf(d+1),1,\smhalf(d+2),1-b,-a\Big).$$
Similar steps lead to
$$\Dsc(a,b)=\frac{\sqrt{b\pi}\,\Gamma\big(\smhalf(d+3)\big)}
{2\Gamma\big(\smhalf(d+4)\big)}\Phi_1\Big(\smhalf(d+3),1,\smhalf(d+4),1-b,-a\Big).$$
Substitution of these alternative $\Csc(a,b)$ and $\Dsc(a,b)$ expressions
into (\ref{eq:goodEgg}) then gives
$$\nabla_{\by}\log\{\pDens(\by)\}=-\by\,
\frac{(d+1)\Phi_1\Big(\smhalf(d+3),1,\smhalf(d+4),1-\tau^2,-\smhalf\Vert\by\Vert^2\Big)}
{(d+2)\Phi_1\Big(\smhalf(d+1),1,\smhalf(d+2),1-\tau^2,-\smhalf\Vert\by\Vert^2\Big)}.
$$
From Appendix B of Gordy (1998), and noting that $0<1-\tau^2<1$,
\begin{equation}
\begin{array}{l}
\Phi_1\Big(\smhalf(d+3),1,\smhalf(d+4),1-\tau^2,-\smhalf\Vert\by\Vert^2\Big)\\[1ex]
\qquad\qquad=\exp\big(-\smhalf\Vert\by\Vert^2\big)
{\displaystyle\sum_{n=0}^{\infty}}
{\displaystyle\frac{\big(\smhalf(d+3)\big)_n
(1-\tau^2)^n}{\big(\smhalf(d+4)\big)_n}}
\oneFone\big(\smhalf, \smhalf(d+4)+n,\smhalf\Vert\by\Vert^2\big).
\end{array}
\label{eq:BabyJune}
\end{equation}
Next note that, from Section 13.1.5 of (\ref{eq:StegunRes}),
$$\oneFone\big(\smhalf, \smhalf(d+4)+n, \smhalf\Vert\by\Vert^2\big)=
\frac{\Gamma\big(\smhalf(d+4)+n\big)}{\Gamma(\smhalf)}
\exp\big(\smhalf\Vert\by\Vert^2\big)\big(\smhalf\Vert\by\Vert^2\big)^{-\smhalf(d+3)-n}
\big\{1+O(\Vert\by\Vert^{-2})\big\}.
$$
Substitution into (\ref{eq:BabyJune}) then gives
{\setlength\arraycolsep{1pt}
\begin{eqnarray*}
&&\Phi_1\Big(\smhalf(d+3),1,\smhalf(d+4),1-\tau^2,-\smhalf\Vert\by\Vert^2\Big)\\[1ex]
&&\quad=\frac{2^{(d+3)/2}\Gamma\big(\smhalf(d+4)\big)}
{\Gamma\big(\smhalf(d+3)\big)\sqrt{\pi}}\Vert\by\Vert^{-(d+3)}\big\{1+O(\Vert\by\Vert^{-2})\big\}
\sum_{n=0}^{\infty}\frac{\Gamma\big(\smhalf(d+3)+n\big)\{2(1-\tau^2)\}^n}{\Vert\by\Vert^{2n}}\\[1ex]
&&\quad=\frac{2^{(d+3)/2}\Gamma\big(\smhalf(d+4)\big)}
{\sqrt{\pi}}\Vert\by\Vert^{-(d+3)}\big\{1+O(\Vert\by\Vert^{-2})\big\}.
\end{eqnarray*}
}
Similarly,
\begin{equation}
\begin{array}{l}
\Phi_1\Big(\smhalf(d+1),1,\smhalf(d+2),1-\tau^2,-\smhalf\Vert\by\Vert^2\Big)\\[1ex]
\qquad\qquad=\exp\big(-\smhalf\Vert\by\Vert^2\big)
{\displaystyle\sum_{n=0}^{\infty}}
{\displaystyle\frac{\big(\smhalf(d+1)\big)_n(1-\tau^2)^n}{\big(\smhalf(d+2)\big)_n}}
\oneFone\big(\smhalf, \smhalf(d+2)+n,\smhalf\Vert\by\Vert^2\big).
\end{array}
\label{eq:BabyJuly}
\end{equation}
Again, from (\ref{eq:StegunRes}),
$$\oneFone\big(\smhalf, \smhalf(d+2)+n,\smhalf\Vert\by\Vert^2\big)=
\frac{\Gamma\big(\smhalf(d+2)+n\big)}{\Gamma(\smhalf)}
\exp\big(\smhalf\Vert\by\Vert^2\big)\big(\smhalf\Vert\by\Vert^2\big)^{-\smhalf(d+1)-n}
\big\{1+O(\Vert\by\Vert^{-2})\big\}.
$$
Substitution into (\ref{eq:BabyJuly}) leads to 
{\setlength\arraycolsep{1pt}
\begin{eqnarray*}
&&\Phi_1\Big(\smhalf(d+1),1,\smhalf(d+2),1-\tau^2,-\smhalf\Vert\by\Vert^2\Big)\\[1ex]
&&\quad=\frac{2^{(d+1)/2}\Gamma\big(\smhalf(d+2)\big)}
{\Gamma\big(\smhalf(d+1)\big)\sqrt{\pi}}\Vert\by\Vert^{-(d+1)}\big\{1+O(\Vert\by\Vert^{-2})\big\}
\sum_{n=0}^{\infty}\frac{\Gamma\big(\smhalf(d+1)+n\big)\{2(1-\tau^2)\}^n}
{\Vert\by\Vert^{2n}}\\[1ex]
&&\quad=\frac{2^{(d+1)/2}\Gamma\big(\smhalf(d+2)\big)}
{\sqrt{\pi}}\Vert\by\Vert^{-(d+1)}\big\{1+O(\Vert\by\Vert^{-2})\big\}.
\end{eqnarray*}
}
We then have
$$
\frac{\Phi_1\Big(\smhalf(d+3),1,\smhalf(d+4),1-\tau^2,-\smhalf\Vert\by\Vert^2\Big)}
{\Phi_1\Big(\smhalf(d+1),1,\smhalf(d+2),1-\tau^2,-\smhalf\Vert\by\Vert^2\Big)}
=
\frac{(d+2)\Vert\by\Vert^{-2}\big\{1+O(\Vert\by\Vert^{-2})\big\}}
{\{1+O(\Vert\by\Vert^{-2})\}}.
$$
This leads to
\begin{equation}
\nabla_{\by}\log\{\pDens(\by)\}=\frac{-(d+1)\by\Vert\by\Vert^{-2}
\big\{1+O\big(\Vert\by\Vert^{-2}\big)\big\}}{1+O\big(\Vert\by\Vert^{-2}\big)}
\quad\mbox{for}\quad \Vert\by\Vert\gg1.
\label{eq:tauLT1}
\end{equation}

\subsection{Tail Limit of the Score Function}

Results (\ref{eq:tauEQ1}), (\ref{eq:tauGT1}) and (\ref{eq:tauLT1}) imply that,
for all $\tau>0$, the score has the following leading term tail behaviour:
\begin{equation}
\nabla_{\by}\log\{\pDens(\by)\}\sim\,-\frac{(d+1)\by}{\Vert\by\Vert^2}
\quad\mbox{for}\quad \Vert\by\Vert\gg1.
\label{eq:MajorsCreekPub}
\end{equation}
It follows immediately that 
$$\lim_{\Vert\by\Vert\to\infty}\nabla_{\by}\log\{\pDens(\by)\}=\bzero.$$

\subsection{Explicit Expressions for $E(\btheta|\by)$}

Arguments similar to those used in Section \ref{sec:BCHFreps} for the score function lead to 
$$E(\btheta|\by)=
\frac{
\Phi_1\Big(\threehalves,1,\smhalf(d+4),1-\tau^{-2},\smhalf\Vert\by\Vert^2\Big)}
{(d+2)\Phi_1\left(\smhalf,1,\smhalf(d+2),1-\tau^{-2},\smhalf\Vert\by\Vert^2\right)}\,\by.
$$
An alternative expression, that uses an integration by parts step as described
in the \textit{Proof of Theorem 2} section of Carvalho \textit{et al.} (2010),
is
\begin{equation}
E(\btheta|\by)=
\left\{1-
\frac{(d+1)
\Phi_1\Big(\smhalf,1,\smhalf(d+4),1-\tau^{-2},\smhalf\Vert\by\Vert^2\Big)}
{(d+2)\Phi_1\left(\smhalf,1,\smhalf(d+2),1-\tau^{-2},\smhalf\Vert\by\Vert^2\right)}
\right\}\,\by
\label{eq:KMcColl}
\end{equation}

\subsection{Bounding of $\big\Vert\by-E(\btheta|\by)\big\Vert$}\label{sec:Bounding}

It follows from (\ref{eq:suxEggs}) and (\ref{eq:KMcColl}) that
\begin{equation}
\big\Vert\by-E(\btheta|\by)\big\Vert=\big\Vert \nabla_{\by}\log\{\pDens(\by)\}\big\Vert.
\label{eq:Vacaville}
\end{equation}
Because of (\ref{eq:MajorsCreekPub}) we then have
$$\nabla_{\by}\log\{\pDens(\by)\}=\bzero\ \mbox{at}\ \by=\bzero\quad\mbox{and}\quad
\nabla_{\by}\log\{\pDens(\by)\}\sim -\frac{(d+1)\by}{\Vert\by\Vert^2}\quad\mbox{for}\quad\Vert\by\Vert\gg1.$$
Relationship (\ref{eq:Vacaville}) then provides
$$\big\Vert\by-E(\btheta|\by)\big\Vert=0\  \mbox{at}\ \by=\bzero\quad\mbox{and}\quad
\big\Vert\by-E(\btheta|\by)\big\Vert\sim\frac{(d+1)}{\Vert\by\Vert}
\quad\mbox{for}\quad\Vert\by\Vert\gg1.
$$
This result and the fact that $\big\Vert\by-E(\btheta|\by)\big\Vert$ is continuous 
in $\by$ implies that $E\big\Vert\by-E(\btheta|\by)\big\Vert$ is bounded by 
some $b_{\tau}<\infty$ that depends only on $\tau$.

\section{Derivation of Result \ref{res:riskRes}}\label{sec:derivResFour}

Consider the Bayesian model 
$$\pDens(\by_1,\ldots,\by_n|\btheta)=\prod_{i=1}^n \pDens(\by_i|\btheta),\quad
\btheta\ \mbox{has prior density function}\ \pHSd\big(\btheta\big),$$
where, for $1\le i\le n$,
$$\pDens(\by_i|\btheta)\equiv(2\pi\sigma^2)^{-d/2}\exp\left\{
-\frac{\Vert\by_i-\btheta\Vert^2}{2\sigma^2}\right\}.$$
Suppose that $\btheta^0$ is the true value of $\btheta$. For each $\btheta\in\real^d$, 
the Kullback-Leibler divergence of $\pDens(\by_i|\btheta)$ from $\pDens(\by_i|\btheta^0)$
is
$$\KL\Big(\pDens(\by_i|\btheta^0)\Vert 
\pDens(\by_i|\btheta) \Big)=\frac{\Vert\btheta-\btheta_0\Vert^2}{2\sigma^2}.$$
For each $\varepsilon>0$ let 
$$A_{\varepsilon}\equiv\left\{\btheta:\KL\Big(\pDens(\by_i|\btheta^0)\Vert 
\pDens(\by_i|\btheta)\Big)\le\varepsilon\right\}
=\{\Vert\btheta-\btheta^0\Vert^2\le 2\sigma^2\varepsilon\}.
$$
Application Proposition 1 of Bhadra \textit{et al.} (2017), which is established in 
Barron (1987), with $\varepsilon=1/n$ leads to
\begin{equation}
R_n\le\frac{1}{n}-\frac{1}{n}\log\left(\int_{A_{1/n}} \pHSd(\btheta)\,d\btheta\right)
=\frac{1}{n}-\frac{1}{n}
\log\left(
\int_{\Vert\btheta-\btheta^0\Vert\le\sigma\sqrt{2/n}} \pHSd(\btheta)\,d\btheta
\right).
\label{eq:CarmenH}
\end{equation}
\subsection{The $\btheta^0=\bzero$ Case}
For the $\btheta^0=\bzero$ case (\ref{eq:CarmenH}) reduces to 
\begin{equation}
R_n\le \frac{1}{n}-\frac{1}{n}
\log\left(\int_{\Vert\btheta\Vert\le\sigma\sqrt{2/n}} \pHSd(\btheta)\,d\btheta\right).
\label{eq:RetaH}
\end{equation}
To determine the order of magnitude of the right-hand side of 
(\ref{eq:RetaH}) we consider separately (a) $d=1$ and (b) $d\ge2$.

\subsubsection{The $d=1$ Case}

When $d=1$ the bound in (\ref{eq:RetaH}) becomes 
\begin{equation}
{\setlength\arraycolsep{1pt}
\begin{array}{rcl}
R_n&\le&\dsoon-\dsoon
\log\left(
\sqrt{2/\pi^3}\,{\displaystyle\int_{0}^{\sigma\sqrt{2/n}}} \exp(\theta^2/2)\,
E_1(\theta^2/2)\,d\theta\right)\\[3ex]
&=&\dsoon-\dsoon\log\left(\pi^{-3/2}
{\displaystyle\int_{0}^{\sigma^2/n}} u^{-1/2}\exp(u)\,
E_1(u)\,du\right).
\end{array}
}
\label{eq:StarOfTheSea}
\end{equation}
From Section 8.214 of Gradshteyn \myand Ryzhik (1994),
\begin{equation}
E_1(u)=-\gamma-\log(u)-\sum_{k=1}^{\infty}\frac{(-u)^k}{k(k!)}
\label{eq:MuffinGuitar}
\end{equation}
where $\gamma\equiv-\mydigamma(1)$ is Euler's constant.
Combining (\ref{eq:MuffinGuitar}) with the Taylor series expansion of $\exp(u)$ 
we obtain
{\setlength\arraycolsep{1pt}
\begin{eqnarray*}
\int_{0}^{\sigma^2/n} u^{-1/2}\exp(u)\,E_1(u)\,du
&=&\int_{0}^{\sigma^2/n}\big\{-\gamma u^{-1/2}-u^{-1/2}\log(u)\\
&&\qquad\quad+(1-\gamma)u^{1/2}
-u^{1/2}\log(u)+\quarter(3-2\gamma)u^{3/2}+\ldots\big\}\,du\\[1ex]
&=&2\sigma\log(n)\,n^{-1/2}+O(n^{-1/2}).
\end{eqnarray*}
}
It follows that 
$$\log\left(\pi^{-3/2}
{\displaystyle\int_{0}^{\sigma^2/n}} u^{-1/2}\exp(u)\,
E_1(u)\,du\right)=-\smhalf\log(n)+\log\{\log(n)\}+O(1).$$
Substitution into (\ref{eq:StarOfTheSea}) then leads to 
$$R_n\le\frac{\log(n)}{2n}-\frac{\log\{\log(n)\}}{n}+O\left(\frac{1}{n}\right).$$

\subsubsection{The $d\ge2$ Case}

In this section we assume that $d\ge2$. To analyze 
$$\gothJ(d,n,\sigma)\equiv\int_{\Vert\btheta\Vert\le\sigma\sqrt{2/n}} \pHSd(\btheta)\,d\btheta$$
we switch to hyper-spherical coordinates as follows:
{\setlength\arraycolsep{1pt}
\begin{eqnarray*}
\theta_1&=&r\cos(\phi_1),\\[1ex]
\theta_2&=&r\sin(\phi_1)\cos(\phi_2),\\[1ex]
\theta_3&=&r\sin(\phi_1)\sin(\phi_2)\cos(\phi_3),\\[1ex]
        &\vdots&     \\[1ex]
\theta_{d-1}&=&r\sin(\phi_1)\cdots\sin(\phi_{d-2})\cos(\phi_{d-1}),\\[1ex]
\theta_d&=&r\sin(\phi_1)\cdots\sin(\phi_{d-2})\sin(\phi_{d-1})\\
\end{eqnarray*}
}
where $r\ge0$, $0\le\phi_1,\phi_2,\ldots,\phi_{d-2}\le\pi$ and $0\le\phi_{d-1}<2\pi$.
Then, noting that $\Vert\btheta\Vert=r$ and the determinant
of the Jacobian of the transformation is such that
$$d\btheta=r^{d-1}\sin^{d-2}(\phi_1)\sin^{d-3}(\phi_2)\cdots\sin(\phi_{d-2})
dr\,d\phi_1\cdots d\phi_{d-1}$$
application of Wallis' Theorem and some additional, but 
straightforward, algebra leads to 
\begin{equation}
\gothJ(d,n,\sigma)=\frac{\Gamma\big((d+1)/2\big)}{\pi\Gamma(d/2)}
\int_0^{\sigma^2/n}u^{-1/2}\exp(u)E_{(d+1)/2}(u)\,du
\label{eq:SloveniaFrontier}
\end{equation}
The next step involves approximation of the integral in 
(\ref{eq:SloveniaFrontier}) using series expansions of 
$E_{(d+1)/2}(u)$ and then applying results such as
\begin{equation}
{\setlength\arraycolsep{0pt}
\begin{array}{rcl}
&&\displaystyle{\int_0^{\sigma^2/n}}u^{-1/2}\,du=2\sigma n^{-1/2},\\
&&\displaystyle{\int_0^{\sigma^2/n}}u^{1/2}\log(u)\,du
=-\twothirds\sigma^3\log(n)n^{-3/2}+\twoninths
\{6\log(\sigma)-2\}\sigma^3n^{-3/2},\\
&&\displaystyle{\int_0^{\sigma^2/n}}u^{1/2}\,du=\twothirds
\sigma^3 n^{-3/2},\\
&&\displaystyle{\int_0^{\sigma^2/n}}u^{3/2}\log(u)\,du=
-\twofifths\sigma^5\log(n)n^{-5/2}+
\twotwentyfifths\{10\log(\sigma)-2\}\sigma^5n^{-5/2}\\
\mbox{and}\
&&\displaystyle{\int_0^{\sigma^2/n}}u^{3/2}\,du=\twofifths\sigma^5n^{-5/2}.
\end{array}
}
\label{eq:ApricotCroissant}
\end{equation}

\vskip2mm
\leftline{$\underline{\mbox{\textit{The $d$ Odd Case}}}$}
\vskip2mm

If $d$ is odd then $\smhalf(d+1)\in\naturalNumbers$ and, from 
from 8.19.8 of Olver (2023),
$$
E_{(d+1)/2}(u)=\frac{(-u)^{(d-1)/2}
\{\mydigamma\big((d+1)/2\big)-\log(u)\}}{\{(d-1)/2\}!}
-\sum_{k=0,k\ne(d-1)/2}^{\infty}\frac{2(-u)^k}{(2k-d+1)k!}.
$$
This expansion, when combined with the Taylor series expansion of $\exp(u)$,
leads to 
$$
{\setlength\arraycolsep{0pt}
\begin{array}{rcl}
&&\displaystyle{\int_{0}^{\sigma^2/n} u^{-1/2}\exp(u)\,E_{(d+1)/2}(u)\,du}\\[1ex]
&&\quad =\displaystyle{\int_{0}^{\sigma^2/n}}
(u^{-1/2}+u^{1/2}+\smhalf u^{3/2}+\sixth u^{5/2}+\ldots)\\[2ex]
&&\qquad\times\left\{\displaystyle{\frac{(-u)^{(d-1)/2}
\{\mydigamma\big((d+1)/2\big)-\log(u)\}}{\{(d-1)/2\}!}}
-\displaystyle{\sum_{k=0,k\ne(d-1)/2}^{\infty}\frac{2(-u)^k}{(2k-d+1)k!}}\right\}\,du\\[6ex]
&&\quad=\displaystyle{\int_{0}^{\sigma^2/n}}
\Bigg\{\displaystyle{\frac{(-1)^{(d-1)/2}u^{(d-2)/2}
\{\mydigamma\big((d+1)/2\big)-\log(u)\}}{\{(d-1)/2\}!}}\\
&&\qquad\qquad\qquad\qquad\qquad\qquad\qquad\qquad\qquad\qquad
-\displaystyle{\sum_{k=0,k\ne(d-1)/2}^{\infty}\frac{2(-1)^ku^{k-1/2}}{(2k-d+1)k!}}
\Bigg\}\,du\\[6ex]
&&\qquad+\displaystyle{\int_{0}^{\sigma^2/n}}
\Bigg\{\displaystyle{\frac{(-1)^{(d-1)/2}u^{d/2}
\{\mydigamma\big((d+1)/2\big)-\log(u)\}}{\{(d-1)/2\}!}}\\
&&\qquad\qquad\qquad\qquad\qquad\qquad\qquad\qquad\qquad\qquad
-\displaystyle{\sum_{k=0,k\ne(d-1)/2}^{\infty}\frac{2(-1)^ku^{k+1/2}}{(2k-d+1)k!}}
\Bigg\}\,du+\ldots.
\end{array}
}
$$
Application of (\ref{eq:ApricotCroissant}) to the early terms in 
these series of integrals reveals that the leading term of the
integral in (\ref{eq:SloveniaFrontier}) is
$$-\frac{2}{(0 - d +1)0!}\int_0^{\sigma^2/n}u^{-1/2}\,du
=\left(\frac{4\sigma}{d-1}\right)n^{-1/2}.
$$
The second term is
$$O\{n^{-3/2}\log(n)\}\quad\mbox{if $d=3$ and}\quad
O(n^{-3/2})\quad\mbox{if $d=5,7,9,\ldots$.}
$$
Therefore
\begin{equation}
\int_{0}^{\sigma^2/n} u^{-1/2}\exp(u)\,E_{(d+1)/2}(u)\,du=
\displaystyle{\left(\frac{4\sigma}{d-1}\right)\,n^{-1/2}
+O\Big(n^{-3/2}\log(n)^{I(d=3)}\Big)}
\label{eq:oddRes}
\end{equation}
for all odd integers $d$ exceeding $1$.

\vskip2mm
\leftline{$\underline{\mbox{\textit{The $d$ Even Case}}}$}
\vskip2mm

If $d$ is even then $\smhalf(d+1)\notin\naturalNumbers$ and, from 
8.19.11 of Olver (2023),
$$
\Gamma\big((1-d)/2\big)^{-1}\exp(u)E_{(d+1)/2}(u)=u^{(d-1)/2}\exp(u)
-\sum_{k=0}^{\infty}\frac{u^k}{\Gamma\big(k+(3-d)/2\big)}.
$$
Hence 
{\setlength\arraycolsep{1pt}
\begin{eqnarray*}
\Gamma\big((1-d)/2\big)^{-1}
\int_{0}^{\sigma^2/n} u^{-1/2}\exp(u)\,E_{(d+1)/2}(u)\,du
&=&\int_0^{\sigma^2/n}u^{(d-2)/2}\big(1+u+\smhalf u^2+\ldots\big)\,du\\[1ex]
&&\quad-2\sum_{k=0}^{\infty}
\frac{\sigma^{2k+1}n^{-(2k+1)/2}}{(2k+1) \Gamma\big(k+(3-d)/2\big)}\\[1ex]
&=&\frac{-2\sigma n^{-1/2}}{\Gamma\big((3-d)/2\big)}+O\big(n^{-\min(d,3)/2}\big)
\end{eqnarray*}
}
and we have 
\begin{equation}
{\setlength\arraycolsep{1pt}
\begin{array}{rcl}
\displaystyle{\int_{0}^{\sigma^2/n} u^{-1/2}\exp(u)\,E_{(d+1)/2}(u)\,du}
&=&\,-\displaystyle{\frac{2\sigma\Gamma\big((1-d)/2\big)n^{-1/2}}{\Gamma\big(1+(1-d)/2\big)}}
+O\big(n^{-\min(d,3)/2}\big)
\\[3ex]
&=&\displaystyle{\left(\frac{4\sigma}{d-1}\right)}\,n^{-1/2}+O\big(n^{-\min(d,3)/2}\big)
\end{array}
}
\label{eq:evenRes}
\end{equation}
for all even positive integers $d$.
\vskip8mm

\noindent
From (\ref{eq:RetaH}), (\ref{eq:oddRes}) and (\ref{eq:evenRes}) we have
$$R_n\le\frac{1}{n}-\frac{1}{n}
\log\left(\frac{4\Gamma\big((d+1)/2\big)\sigma n^{-1/2}}
{\pi\Gamma(d/2)(d-1)}+O(n^{-1})\right)
=\frac{\log(n)}{2n}+O\left(\frac{1}{n}\right)
$$
for all $d\ge2$.

\subsection{The $\btheta^0\ne\bzero$ Case}\label{sec:thetaNonZeroCase}

In the $\btheta^0\ne\bzero$ case we have the bound 
\begin{equation}
R_n\le \frac{1}{n}-\frac{1}{n}
\log\left(\int_{S(\btheta^0,\sigma,n)} \pHSd(\btheta)\,d\btheta\right).
\label{eq:RetaI}
\end{equation}
where
$$S(\btheta^0,\sigma,n)\equiv\{\btheta: \Vert\btheta-\btheta^0\Vert\le\sigma\sqrt{2/n}\}.$$
If $C(\btheta^0,\sigma,n)$ is the largest hypercube inscribed in $S(\btheta^0,\sigma,n)$ then
\begin{equation}
\int_{S(\btheta^0,\sigma,n)} \pHSd(\btheta)\,d\btheta\ge 
\int_{C(\btheta^0,\sigma,n)} \pHSd(\btheta)\,d\btheta.
\label{eq:Butterscotch}
\end{equation}
For sufficiently large $n$, $\pHSd(\btheta)$ is a very smooth function over $C(\btheta^0,\sigma,n)$
and Taylor series arguments can be used to establish that
\begin{equation}
\int_{C(\btheta^0,\sigma,n)} \pHSd(\btheta)\,d\btheta=K(d,\btheta^0,\sigma)n^{-d/2}\{1+o(1)\}
\label{eq:Castle}
\end{equation}
for some positive constant $K(d,\btheta^0,\sigma)$. Application of (\ref{eq:Butterscotch})
and (\ref{eq:Castle}) to the bound in (\ref{eq:RetaI}) then leads to
$$R_n\le\frac{d\log(n)}{2n}+O\left(\frac{1}{n}\right)\quad\mbox{for}\quad\btheta^0\ne\bzero.$$

\section{Derivation of Result \ref{res:threshSide}}\label{sec:derivResFive}

We first obtain an explicit expression for the
posterior density function of $\lambda$.
Note that
$$\pDens(\lambda|\by)=\frac{\pDens(\by,\lambda)}{\pDens(\by)}
=\frac{\int_{\real^d}\pDens(\by,\bpsi,\lambda)d\bpsi}{\pDens(\by)}
=\frac{\pDens(\lambda)
\int_{\real^d}\pDens(\by|\bpsi)\pDens(\bpsi|\lambda)d\bpsi}{\pDens(\by)}.
$$
It is easy to establish that
$$
\pDens(\by|\bpsi)
=(2\pi\fixOne)^{-d/2}
\exp\left(-\frac{\Vert\by\Vert^2}{2\fixOne}\right)
\exp\left\{(1/\fixOne)
\left[
\begin{array}{c}
\bpsi\\
\vecof(\bpsi\bpsi^T)
\end{array}
\right]^T
\left[
\begin{array}{c}
\by\\
-\smhalf\vecof(\bI_d)
\end{array}
\right]
\right\}
$$
and
$$
\pDens(\bpsi|\lambda)
=(2\pi\lambda^2\fixTwo)^{-d/2}
\exp\left[\{1/(\lambda^2\fixTwo)\}
\left[
\begin{array}{c}
\bpsi\\
\vecof(\bpsi\bpsi^T)
\end{array}
\right]^T
\left[
\begin{array}{c}
\bzero\\
-\smhalf\vecof(\bI_d)
\end{array}
\right]
\right].
$$
Hence 
$$
\pDens(\by|\bpsi)\pDens(\bpsi|\lambda)=(2\pi\fixOne)^{-d/2}
\exp\left(-\frac{\Vert\by\Vert^2}{2\fixOne}\right)
(2\pi\lambda^2\fixTwo)^{-d/2}
\exp\left\{
\left[
\begin{array}{c}
\bpsi\\
\vecof(\bpsi\bpsi^T)
\end{array}
\right]^T
\left[
\begin{array}{c}
\bdeta_1\\
\bdeta_2
\end{array}
\right]
\right\}
$$
where
\begin{equation}
\bdeta_1\equiv\by\big/\fixOne\quad\mbox{and}\quad
\bdeta_2\equiv -\smhalf\left(\frac{1}{\fixOne}+\frac{1}{\lambda^2\fixTwo}\right)
\vecof(\bI_d).
\label{eq:Bacall}
\end{equation}
As a function of $\bpsi$, $\pDens(\by|\bpsi)\pDens(\bpsi|\lambda)$ is
a Multivariate Normal density function with natural
parameters given by (\ref{eq:Bacall}). Therefore, standard results
concerning the normalizing factor of the Multivariate
Normal family leads to 
{\setlength\arraycolsep{1pt}
\begin{eqnarray*}
\int_{\real^d}\pDens(\by|\bpsi)\pDens(\bpsi|\lambda)d\bpsi
&=&(2\pi\fixOne)^{-d/2}\exp\left(-\frac{\Vert\by\Vert^2}{2\fixOne}\right)
(\lambda^2\fixTwo)^{-d/2}
\left|\left(\frac{1}{\fixOne}+\frac{1}{\lambda^2\fixTwo}\right)\bI_d\right|^{-1/2}\\
&&\qquad\times\exp\left[\frac{1}{2(\fixOne)^2}\by^T\left\{
\left(\frac{1}{\fixOne}+\frac{1}{\lambda^2\fixTwo}\right)\bI_d\right\}^{-1}\by\right]\\[1ex]
&=&(2\pi\fixOne)^{-d/2}\left(1+\frac{\lambda^2\fixTwo}{\fixOne}\right)^{-d/2}
\exp\left\{-\frac{\Vert\by\Vert^2}{2(\fixOne+\lambda^2\fixTwo)}\right\}.
\end{eqnarray*}
}
We then have the following expression for the posterior density function of $\lambda$:
$$\pDens(\lambda|\by)=\frac{2I(\lambda>0)}{\pi(1+\lambda^2)\pDens(\by)}
\left\{2\pi(\fixOne+\lambda^2\fixTwo)\right\}^{-d/2}
\exp\left\{-\frac{\Vert\by\Vert^2}{2(\fixOne+\lambda^2\fixTwo)}\right\}.
$$

The posterior density function of $\bpsi$ is 
$$\pDens(\bpsi|\by)=\frac{\pDens(\by,\bpsi)}{\pDens(\by)}
=\frac{\int_0^{\infty}\pDens(\by,\bpsi,\lambda)d\lambda}{\pDens(\by)}
=\frac{\pDens(\by|\bpsi)\int_0^{\infty}\pDens(\bpsi|\lambda)
\pDens(\lambda)d\lambda}{\pDens(\by)}.
$$
Introduce following the notation, defined immediately after equation (4) 
of Wand \myand Jones (1993):
$$\phi_{\bSigma}(\bx)=(2\pi)^{-d/2}|\bSigma|^{-1/2}\exp(-\smhalf\bx^T\bSigma^{-1}\bx).$$
Then
{\setlength\arraycolsep{1pt}
\begin{eqnarray*}
\int_0^{\infty}\pDens(\bpsi|\lambda)
\pDens(\lambda)d\lambda
&=&\int_0^{\infty}(2\pi\lambda^2\fixTwo)^{-d/2}
\exp\left(-\frac{\Vert\bpsi\Vert^2}{2\lambda^2\fixTwo}\right)
\frac{2}{\pi(1+\lambda^2)}\,d\lambda\\[1ex]
&=&(2/\pi)\int_0^{\infty}\phi_{\lambda^2\fixTwo\bI_d}(\bpsi)
\frac{1}{(1+\lambda^2)}\,d\lambda.
\end{eqnarray*}
}
Since
$$\pDens(\by|\bpsi)=\phi_{\fixOne\bI_d}(\by-\bpsi)=\phi_{\fixOne\bI_d}(\bpsi-\by)$$
we then have
\begin{equation}
\pDens(\by)\pDens(\bpsi|\by)=(2/\pi)\int_0^{\infty}
\phi_{\fixOne\bI_d}(\bpsi-\by)
\phi_{\lambda^2\fixTwo\bI_d}(\bpsi-\bzero)
\frac{1}{(1+\lambda^2)}\,d\lambda.
\label{eq:Zuppe}
\end{equation}
From (A.1) of Wand \myand Jones (1993),
{\setlength\arraycolsep{1pt}
\begin{eqnarray*}
\phi_{\fixOne\bI_d}(\bpsi-\by)
\phi_{\lambda^2\fixTwo\bI_d}(\bpsi-\bzero)
&=&\phi_{(\fixOne+\lambda^2\fixTwo)\bI_d}(\by)
\phi_{\{\lambda^2\fixOne\fixTwo/(\fixOne+\lambda^2\fixTwo)\}\bI_d}
\left(\bpsi-\frac{\lambda^2\fixTwo\by}{\fixOne+\lambda^2\fixTwo}\right)\\[1ex]
&=&\{2\pi(\fixOne+\lambda^2\fixTwo)\}^{-d/2}
\exp\left\{-\frac{\Vert\by\Vert^2}{2(\fixOne+\lambda^2\fixTwo)}\right\}\\[1ex]
&&\qquad\times\phi_{\{\lambda^2\fixOne\fixTwo/(\fixOne+\lambda^2\fixTwo)\}\bI_d}
\left(\bpsi-\frac{\lambda^2\fixTwo\,\by}{\fixOne+\lambda^2\fixTwo}\right).
\end{eqnarray*}
}
Therefore, the posterior density function of $\bpsi$ has the following expression
in terms of the posterior density function of $\lambda$:
$$
\pDens(\bpsi|\by)=\int_0^{\infty}\pDens(\lambda|\by)
\phi_{\{\lambda^2\fixOne\fixTwo/(\fixOne+\lambda^2\fixTwo)\}\bI_d}
\left(\bpsi-\frac{\lambda^2\fixTwo\,\by}{\fixOne+\lambda^2\fixTwo}\right)\,d\lambda.
$$
Hence, with an interchange in the order of integration,
{\setlength\arraycolsep{1pt}
\begin{eqnarray*}
E(\bpsi|\by)&=&\int_0^{\infty}\pDens(\lambda|\by)
\left\{\int_{\real^d}\bpsi
\phi_{\{\lambda^2\fixOne\fixTwo/(\fixOne+\lambda^2\fixTwo)\}\bI_d}
\left(\bpsi-\frac{\lambda^2\fixTwo\,\by}{\fixOne+\lambda^2\fixTwo}\right)d\bpsi\right\}\,d\lambda\\[1ex]
&=&\int_0^{\infty}\pDens(\lambda|\by)\left(\frac{\lambda^2\fixTwo}
{\fixOne+\lambda^2\fixTwo}\right)\,d\lambda\,\by
=E\left(\frac{\lambda^2\,\fixTwo}{\fixOne+\lambda^2\fixTwo}\Big|\by\right)\by
\end{eqnarray*}
}
as required.

\section*{References}

\bib
Abramowitz, M. \myand Stegun, I.A. eds. (1968).
\textit{Handbook of Mathematical Functions}.
New York: Dover.

\bib
Barron, A.R. (1987). Are Bayes rules consistent in information?
In T.M. Cover \myand B. Gopinath (eds.),
\textit{Open Problems in Communication and Computation}, pp. 85--91,
New York: Springer-Verlag.

\bib
Bhadra, A., Datta, J., Polson, N.G. \myand Willard, B. (2017).
The horseshoe+ estimator for ultra-sparse signals.
\textit{Bayesian Analysis}, \textbf{12}, 1105--1131.

\bib
Gordy, M.B. (1998). A generalization of generalized beta distributions.
In \textit{Finance and Economics Discussion Series.} Board of Governors
of the Federal Reserve System, United States of America.

\bib
Gradshteyn, I.S. \myand Ryzhik, I.M. (1994). 
\textit{Tables of Integrals, Series, and Products}, 5th Edition,
San Diego, California: Academic Press.

\bib
Olver, F.W., Olde Daalhuis, D., Lozier, D.W., 
Schneider, B.I., Boisvert, C.W., Clark, C.W., Miller, B.R.,
Saunders, B.V., Cohl, H.S. \myand McClain, M.A. (eds.). (2023). 
\textit{U.S. National Institute of Standards and Technology
Digital Library of Mathematical Functions.}
Release 1.1.11.\\
\texttt{https://dlmf.nist.gov}

\bib
Wand, M.P. and Jones, M.C. (1993). Comparison of smoothing
parameterizations in bivariate density estimation.
{\it Journal of the American Statistical
Association}, {\bf 88}, 520--528.

\end{document}